\documentclass{article}
\usepackage[]{epsf,epsfig,amsmath,amssymb,amsfonts,latexsym}

\newtheorem{thma}{Theorem}[section]
\newtheorem{theorem}{Theorem}[section]
\newtheorem{lemma}[theorem]{Lemma}

\newtheorem{fact}[theorem]{Fact}

\newtheorem{conjecture}[theorem]{Conjecture}
\newtheorem{proposition}[theorem]{Proposition}

\newtheorem{corollary}[theorem]{Corollary}

\newtheorem{remark}[theorem]{Remark}

\def\eps{\varepsilon}

\def \PP {\mathbb P}
\def\qed{\hfill $\vcenter{\hrule height .3mm
\hbox {\vrule width .3mm height 2.1mm \kern 2mm \vrule width .3mm
height 2.1mm} \hrule height .3mm}$ \bigskip}
\def\P{\mathbb{P}}

\def\EE{\mathbb{E}}

\def\RR{\mathbb{R}}

\def\Sph{S^{n-1}}

\begin{document}
\title{Thin Shell Implies Spectral Gap up to Polylog via a Stochastic Localization Scheme}
\author{Ronen Eldan\thanks{Supported in part by the Israel Science Foundation and by a Marie Curie Grant from the Commission of the European Communities.}}

\maketitle
\abstract{We consider the isoperimetric inequality on the class of high-dimensional isotropic convex bodies. We establish quantitative connections between two well-known open problems related to this inequality, namely, the \emph{thin shell} conjecture, and the conjecture by \emph{Kannan, Lov\'{a}sz, and Simonovits}, showing that the corresponding optimal bounds are equivalent up to logarithmic factors. In particular we prove that, up to logarithmic factors, the minimal possible ratio between surface area and volume is attained on ellipsoids. We also show that a positive answer to the thin shell conjecture would imply an optimal dependence on the dimension in a certain formulation of the Brunn-Minkowski inequality. Our results rely on the construction of a stochastic localization scheme for log-concave measures. }

\section{Introduction}
The starting point of this paper is a conjecture by  Kannan, Lov\'{a}sz, and Simonovits (in short, the KLS conjecture) about the isoperimetric inequality for convex bodies in $\RR^n$. Roughly speaking, The KLS conjecture asserts that, up to a universal constant, the most efficient way to cut a convex body into two parts is with a hyperplane. To be more precise, given convex body $K \subset \RR^n$ whose
barycenter is at the origin, and
a subset $T \subset K$ with $Vol_n(T) = R Vol_n(K)$, the KLS conjecture suggests that
\begin{equation} \label{eqisop}
Vol_{n-1}(\partial T \cap Int(K)) \geq R C \inf_{\theta \in \Sph} Vol_{n-1} (K \cap \theta^\perp)
\end{equation}
for some universal constant $C>0$, whenever $R \leq \frac{1}{2}$. Here, $Vol_{n-1}$ stands for the $(n-1)$-dimensional volume, $\Sph$ is the unit sphere, $\theta^\perp$ is the hyperplane passing through the origin whose normal direction is $\theta$ and $Int(K)$ is the interior of $K$. \\ \\
The point of this paper is to reduce this conjecture to the case where $T$ is an ellipsoid, up to a logarithmic correction. \\ \\
In order to give a precise formulation of the KLS conjecture, we begin with some notation. A probability density $\rho: \RR^n \rightarrow [0, \infty)$
is called {\it log-concave} if it takes the form $\rho = \exp(-H)$
for a convex function $H: \RR^n \rightarrow \RR \cup \{ \infty \}$.
A probability measure is log-concave if it has a log-concave density.
The uniform probability measure on a convex body is an example for a
log-concave probability measure, as well as, say, the gaussian
measure in $\RR^n$. A log-concave probability density decays
exponentially at infinity, and thus has moments of all orders. For a
probability measure $\mu$ on $\RR^n$ with finite second moments, we
consider its barycenter $b(\mu) \in \RR^n$ and covariance matrix $Cov(\mu)$
defined by
$$ b(\mu) = \int_{\RR^n} x d \mu(x), \ \ \ \ \ \ Cov(\mu) = \int_{\RR^n} (x - b(\mu)) \otimes (x - b(\mu)) d \mu(x) $$
where for $x \in \RR^n$ we write $x \otimes x$ for the $n\times n$
matrix $(x_i x_j)_{i,j=1,\ldots,n}$. A log-concave probability
measure $\mu$ on $\RR^n$ is {\it isotropic} if its barycenter lies
at the origin and its covariance matrix is the identity matrix. \\

Given a measure $\mu$, Minkowski's boundary measure of a Borel set $A \subset \RR^n$,
is defined by,
$$
\mu^+(A) = \liminf_{\eps \to 0^+} \frac{\mu(A_\eps) - \mu(A)}{\eps}
$$
where
$$
A_\eps := \{x \in \RR^n; ~~ \exists y, ~ |x-y| \leq \eps \}
$$
is the $\eps$-extension of $A$. \\ \\
The main point of this paper is to find an upper bound for the constant,
\begin{equation} \label{defg}
G_n^{-1} := \inf_{\mu} \inf_{A \subset \RR^n} \frac{\mu^+(A)}{\mu(A)}
\end{equation}
where $\mu$ runs over all isotropic log-concave measures in $\RR^n$ and $A \subset \RR^n$
runs over all Borel sets with $\mu(A) \leq \frac{1}{2}$. \\

The constant $G_n$ is known as the optimal inverse \emph{Cheeger} constant. According to a result of Ledoux, \cite{ledoux}, the quantity $G_n^{-2}$ is also equivalent, up to a universal constant, to the optimal \emph{spectral gap} constant of isotropic log-concave measures in $\RR^n$ (see (\ref{gnpoincare}) below). For an extensive review of this constant and equivalent formulations, see \cite{Mil}. One property of $G_n$ of particular importance in this note is,
\begin{equation} \label{gnpoincare}
\frac{1}{C} G_n^2 \leq \sup_{\mu} \sup_{\varphi} \frac{\int \varphi^2 d \mu}{ \int |\nabla \varphi|^2 d \mu} \leq C G_n^2
\end{equation}
Where $\mu$ runs over all isotropic log-concave measures and $\varphi$ runs over all smooth enough functions with $\int \varphi d \mu = 0$ and
$C>0$ is some universal constant.

In \cite{KLS}, it is conjectured that,
\begin{conjecture}
There exists a universal constant $C$ such that $G_n < C$ for all $n \in \mathbb{N}$.
\end{conjecture}
In this note we will show that, up to a small correction, the above is implied by a seemingly weaker hypothesis. \\

Next, we would like to formulate the \emph{thin-shell conjecture}. Let $\sigma_n \geq 0$ satisfy
\begin{equation}
 \sigma_n^2 = \sup_X  \EE \left [(|X| - \sqrt{n})^2 \right ] \label{eq_1356}
\end{equation}
where the supremum runs over all isotropic, log-concave random vectors $X$ in $\RR^n$. 
The \textbf{shin-shell conjecture} (see Anttila, Ball and Perissinaki \cite{ABP} and Bobkov and Koldobsky \cite{BK}) asserts the following:
\begin{conjecture}
There exists a universal constant $C$ such that, 
\begin{equation} \label{thinshell}
\sigma_n < C
\end{equation}
for all $n \in \mathbb{N}$.
\end{conjecture}

An application of (\ref{gnpoincare}) with the function $\varphi(x) = |x|^2$ shows
that the thin-shell conjecture is \emph{weaker} than the KLS conjecture. \\

The first nontrivial bound for $\sigma_n$ was given by Klartag in \cite{K1}, who showed that $\sigma_n \leq C \frac{n^{1/2}}{\log (n+1)}$. Several improvements have been introduced around the same method, see e.g \cite{K2} and \cite{fleury}. The best known bound for $\sigma_n$ at the time of this note is due to Guedon and E. Milman, in \cite{GM}, extending previous works of Klartag, Fleury and Paouris, who show that $\sigma_n \leq C n^{\frac{1}{3}}$. The thin-shell conjecture was shown to be true for several specific classes of convex bodies,
such as bodies with a symmetry for coordinate reflections (Klartag, \cite{K_uncond}) and certain \emph{random} bodies (Fleury, \cite{fleury2}). \\

It was found by Sudakov, \cite{sudakov}, that the parameter $\sigma_n$ is highly related to almost-gaussian behaviour
of certain marginals of a convex body, a fact now known as the central limit theorem for convex sets \cite{K1}.
This theorem asserts that most of the one-dimensional marginals of an isotropic, 
 log-concave random vector are approximately gaussian in the sense that the Kolmogorov distance to the standard
gaussian distribution of a typical marginal has roughly
the order of magnitude of $\sigma_n / \sqrt{n}$. Therefore the conjectured 
bound (\ref{thinshell}) actually concerns the quality of
the gaussian approximation to the marginals of high-dimensional
log-concave measures.
\\ \\
The first theorem of this note reads, \\
\begin{thma} \label{mainthm1}
There exists a constant $C>0$ such that for all $n \geq 2$,
$$
G_n \leq C \sqrt{(\log n) \sum_{k=1}^n \frac{\sigma_k^2}{k}}.
$$
\end{thma}
Note that, in particular, for any constant $\kappa > 0$ such that $\sigma_n \leq n^{\kappa}$ for all $n \in \mathbb{N}$, one gets $G_n \leq C (\sqrt{\log n}) n^{\kappa}$. Under the \emph{thin-shell conjecture}, the theorem gives $G_n < C \log n$.
\bigskip
\begin{remark}
Plugging the results of this paper into the currently best known bound for $\sigma_n$ (proven in \cite{GM}), $\sigma_n \leq C n^{1/3}$, it follows that
$$
G_n \leq C n^{1/3} \sqrt{\log n}.
$$
This slightly improves the previous bound, $G_n \leq C n^{5/12}$, which is a corollary of \cite{GM} and \cite{Bobkov}.
\end{remark}
\begin{remark}
In \cite{EK1}, B. Klartag and the author have found a connection between the thin-shell hypothesis and another well known conjecture
related to convex bodies, known as the \emph{hyperplane conjecture}. The methods of this paper share some common lines with the methods
in \cite{EK1}. In a very recent paper of K.Ball and V.H. Nguyen, \cite{bn}, a connection between the KLS conjecture and the hyperplane conjecture that applies for individual log-concave measures has also been established. They show that the isotropic constant of a log concave measure which attains a spectral gap is bounded by a constant which depends exponentially on the spectral gap.
\end{remark}
\begin{remark}
Compare this result with the result in \cite{Bobkov}. Bobkov's theorem states that for any log-concave random vector $X$ and any smooth function $\varphi$, one has
$$
\frac{Var[ \varphi(X)]}{ \EE \left [ |\nabla \varphi(X)|^2 \right] } \leq C \EE[|X|] \sqrt{ Var[|X|]}.
$$
Under the \emph{thin-shell} hypothesis, Bobkov's theorem gives $G_n \leq C n^{1/4}$.
\end{remark}
\bigskip
The bound in theorem \ref{mainthm1} will rely on the following intermediate constant which corresponds to a slightly stronger \emph{thin shell} bound. Define,
\begin{equation} \label{defkn}
K_n^2 := \sup_X \sup_{\theta \in \Sph} \sum_{i,j=1}^n \EE[X_i X_j \langle X, \theta \rangle]^2,
\end{equation}
where the supremum runs over all isotropic log-concave random vectors $X$ in $\RR^n$. Obviously, an equivalent definition of $K_n$ will be,
$$
K_n := \sup_{\mu} \left | \left |\int_{\RR^n}x_1 x \otimes x d \mu(x) \right | \right |_{HS}
$$
where the supremum runs over all isotropic log-concave measures in $\RR^n$. Here,
$|| \cdot ||_{HS}$ stands for the Hilbert-Schmidt norm of a matrix. \\ \\
There is a simple relation between $K_n$ and $\sigma_n$, namely,
\begin{lemma} \label{ksigma}
There exists a constant $C>0$ such that for all $n \geq 2$,
$$
K_n \leq C \sqrt{ \sum_{k=1}^n \frac{\sigma_k^2}{k}}.
$$
\end{lemma}
\bigskip
Theorem \ref{mainthm1} will be a consequence of the above lemma along with,
\begin{proposition} \label{gk}
There exists a constant $C>0$ such that for all $n \geq 2$,
$$
G_n \leq C K_n \sqrt{\log n}.
$$
\end{proposition}
\begin{remark}
The constant $K_n$ satisfies the following bound:
$$
K_n^{-1} \geq c \inf_{\mu} \inf_{E \subset \RR^n} \frac{\mu^+(E)}{\mu(E)}
$$
where $\mu$ runs over all isotropic log-concave measures in $\RR^n$, $E$
runs over all \emph{ellipsoids} with $\mu(E) \leq \frac{1}{2}$ and $c>0$ is some universal constant.
This shows that up to the extra factor $\sqrt{ \log n}$, in order to control the
minimal possible surface area among \emph{all possible subsets} of measure $\frac{1}{2}$ on the class of isotropic log-concave measures, it is enough to control the surface area of \emph{ellipsoids}. See section 6 below for details.
\end{remark}

\bigskip
We move on to the second result of this paper, a stability result for the \emph{Brunn-Minkowski Inequality}.
The Brunn-Minkowski inequality states, in one of its normalizations,
that
\begin{equation}
 Vol_n \left( \frac{K + T}{2} \right) \geq \sqrt{Vol_n(K) Vol_n(T)}
 \label{eq_928}
\end{equation}
for any compact sets $K, T \subset \RR^n$, where $(K + T) / 2 = \{
(x + y) / 2 ; x \in K, y \in T \}$ is half of the Minkowski sum of
$K$ and $T$. When $K$ and $T$ are closed convex sets, equality in (\ref{eq_928}) holds if and only if $K$ is a translate of $T$.

\medskip When there is an almost-equality in (\ref{eq_928}), $K$ and $T$ are almost
translates of each other in a certain sense (which varies between different estimates). Estimates of this form,
often referred to as \emph{stability estimates}, appear in Diskant
\cite{diskant}, in Groemer \cite{groemer}, and in Figalli, Maggi and
Pratelli \cite{FMP1, FMP2}, Segal \cite{segal}.

The result \cite{FMP2}, which is essentially the strongest result in its category, and other existing stability estimates
share a common thing: the bounds become worse as the dimension increases. In a recent paper, \cite{EK2}, Klartag
and the author suggested that the correct bounds might actually become better as the dimension increases, as demonstrated
by certain results. The estimates presented here may be viewed as a continuation of this line of research.

In order to formulate our result, we define the two constants
\begin{equation}
\kappa = \liminf_{n \to \infty} \frac{\log \sigma_n}{\log n}, ~~~ \tau_n = \max\left (1, \max_{1 \leq j \leq n} \frac{\sigma_j}{j^\kappa} \right ),
\end{equation}
so that $\sigma_n \leq \tau_n n^\kappa$. Note that the thin-shell conjecture implies $\kappa = 0$ and $\tau_n < C$.  \\ \\
Our main estimate reads,
\begin{thma} \label{mainthm2}
For every $\epsilon > 0$ there exists a constant $C(\epsilon)$ such that the following holds: Let $K,T$ be convex bodies whose volume is $1$ and whose barycenters lie at the origin. Suppose that the covariance matrix of the uniform measure on $K$ is equal to $L_K Id$ for a constant $L_K > 0$. Denote,
\begin{equation} \label{defV1}
V = Vol_n \left ( \frac{K+T}{2} \right ),
\end{equation}
and define
$$
\delta = C(\epsilon) L_K V^5 \tau_n n^{2 (\kappa - \kappa^2) + \epsilon}.
$$
Then,
$$
Vol_n( K_\delta \cap T ) \geq 1 - \epsilon.
$$
\end{thma}
\bigskip
Some remarks:
\begin{remark}
It follows from theorem 1.4 in \cite{EK2} that the above estimate is true with 
$$
\delta =  C(\epsilon) \sqrt{\tau_n} n^{1/4 + \kappa / 2} V^{5/2}.
$$
If $\kappa \geq 1/4$, then the result we prove here weaker than the one in \cite{EK2}. However, under the thin shell hypothesis, the result of this paper becomes stronger, and is in fact tight up to the term $C(\epsilon) n^{\epsilon}$. This tightness is demonstrated, for instance, by taking $K$ and $T$ to be the unit cube and a unit cube truncated by a ball of radius $\sqrt n$ and normalized to be isotropic.
\end{remark}
\begin{remark}
Using the bound in \cite{GM}, the theorem gives 
$$
\delta = C(\epsilon) n^{\frac{4}{9} + \epsilon} V^5 L_K.
$$ 
Note that if the assumption (\ref{defV1}) is dropped, even if the covariance matrices of $K$ and $T$ are assumed to be equal, the best corresponding bound would be $\delta = C \sqrt n L_K$ as demonstrated, for example, by a cube and a ball.
\end{remark}
\begin{remark}
The above bound complements, in some sense, the result proven in \cite{FMP1}, which reads,
$$
Vol_n((K + x_0) \Delta T)^2 \leq n^7 ( Vol_n((K+T)/2) - 1)
$$
for some choice of $x_0$, where $\Delta$ denotes the symmetric difference between the sets.
Unlike the result presented in this paper, the result in \cite{FMP1} gives much more information as the expression $Vol_n((K+T)/2) - 1$ 
approaches zero. On the other hand the result presented here already gives some information when $Vol_n((K+T)/2) = 10$.
\end{remark}
\bigskip
The structure of this paper is as follows: In section 2, we construct a stochastic localization scheme which will be the main ingredient our proofs. In section 3, we establish a bound for the covariance matrix of the measure throughout the localization process, which will be essential for its applications. In section 4, we prove theorem \ref{mainthm1} and in section 5 we prove theorem \ref{mainthm2} and its corollaries. In section 6 we tie some loose ends. \\

\medskip Throughout this note, we use the letters $c,
\tilde{c}, c^{\prime}, C, \tilde{C}, C^{\prime}, C''$ to denote positive
universal constants, whose value is not necessarily the same in
different appearances. Further notation used throughout
the text: for a Borel measure $\mu$ on $\RR^n$, $supp(\mu)$ is the
minimal closed set of full measure. The Euclidean unit ball is denoted by $B_n = \{ x
\in \RR^n ; |x| \leq 1 \}$. Its boundary is denoted by $\Sph$. We write $\nabla \varphi$ for the gradient of the function $\varphi$, and $\nabla^2 \varphi$ for the Hessian matrix. For a positive semi-definite symmetric matrix $A$, we denote its largest eigenvalue by $||A||_{OP}$. For any matrix $A$, we denote the sum of its diagonal entries by $Tr(A)$, and by $||A||_{HS}^2$ we denote the sum of the eigenvalues of the matrix $A^T A$. For two densities $f$, $g$ on $\RR^n$, define the \emph{Wasserstein distance}, $W_2(f, g)$, by
$$
W_2(f, g)^2 = \inf_\xi \int_{\RR^n \times \RR^n} |x - y|^2  d \xi(x,y)
$$
where the infimum is taken over all measures $\xi$ on $\RR^{2n}$ whose marginals onto the first and last $n$ coordinates are the measures whose densities are $f$ and $g$ respectively (see, e.g. \cite{villani} for more information). \\
Finally, for a continuous time stochastic process $X_t$, we denote by $d X_t$ the differential of $X_t$, and by $[X]_t$ the quadratic
variation of $X_t$. For a pair of continuous time stochastic processes $X_t, Y_t$, the quadratic covariation will be denoted by $[X, Y]_t$.
\\ \\

\emph{Acknowledgements}
I owe this work to countless useful discussions I have had with my supervisor, Bo'az Klartag, through which I learnt the vast part of what I know about the subject, as well as about related topics, and for which I am grateful. I would also like to thank Vitali and Emanuel Milman and Boris Tsirelson for inspiring discussions and for their useful remarks on a preliminary version of this note. Finally, I would like to thank the
anonymous referee for doing a tremendous job reviewing a preliminary version of this paper, thanks to his/her ideas the proofs are significantly simpler, shorter and more comprehensible.

\section{A stochastic localization scheme}
In this section we construct the localization scheme which will be the principal component in our proofs.
The construction will use elementary properties of semimartingales and stochastic integration.
For definitions, see \cite{durett}. \\ \\
For the construction, we assume that we are given some isotropic random vector $X \in \RR^n$ with density $f(x)$.
Well-known concentration bounds for log-concave measures (see, e.g., section 2 of \cite{K2}) will allow us to assume throughout the paper that
\begin{equation} \label{compactsupport}
supp(f) \subseteq n B_n,
\end{equation}
where $B_n$ is the Euclidean ball of radius 1. \\ \\
We begin with some definitions. For a vector $c \in \RR^n$ and an $n \times n$ matrix $B$, we write
$$
V_f(c,B) = \int_{\RR^n} e^{\langle c, x \rangle - \frac{1}{2} \langle B x, x \rangle }  f(x) dx.
$$
Define a vector valued function,
$$
a_f(c,B) = V_f^{-1} (c,B) \int_{\RR^n} x e^{\langle c, x \rangle - \frac{1}{2} \langle B x, x \rangle  }  f(x) dx,
$$
and a matrix valued function,
$$
A_f(c,B) = V_f^{-1} (c,B) \int_{\RR^n} (x - a_f(c,B)) \otimes (x - a_f(c,B)) e^{\langle c, x \rangle - \frac{1}{2} \langle B x, x \rangle }  f(x) dx.
$$
The assumption (\ref{compactsupport}) ensures that $V_f$, $a_f$ and $A_f$ are smooth functions of $c, B$. \\ \\
Let $W_t$ be a standard Wiener process and consider the following system of stochastic differential equations:
\begin{equation} \label{stochastic1}
c_0 = 0, ~~ d c_t = A_f^{-1/2}(c_t, B_t) dW_t + A_f^{-1}(c_t, B_t) a_f(c_t, B_t) dt,
\end{equation}
$$
B_0 = 0, ~~ d B_t = A_f^{-1}(c_t, B_t) dt.
$$
Taking into account the fact that the functions $A_f, a_f$ are smooth and that $A_f(c, B)$ is positive definite for all $c,B$,
we can use a standard existence and uniqueness theorem (see e.g., \cite{oksendal}, section 5.2) to ensure the existence and uniqueness of a solution 
in some interval $0 \leq t \leq t_0$, where $t_0$ is an almost-surely positive random variable. \\ \\
Next, we construct a 1-parameter family of functions $\Gamma_t(f)$ by defining, \\
\begin{equation} \label{defineF}
F_t(x) = V_f^{-1} (c_t, B_t) e^{\langle c_t, x \rangle - \frac{1}{2} \langle B_t x, x \rangle }
\end{equation}
and
$$
\Gamma_t(f)(x) = f(x) F_t(x).
$$
Also, abbreviate
$$
a_t = a_f(c_t, B_t), ~~ A_t = A_f(c_t, B_t), ~~ V_t = V_f(c_t, B_t), ~~ f_t = \Gamma_t(f),
$$
so that $a_t$ and $A_t$ are the barycenter and the covariance matrix of the function $f_t$. \\

The following lemma may shed some light on this construction.
\begin{lemma} \label{basicform}
The function $F_t$ satisfies the following set of equations:
\begin{equation} \label{contloc}
F_0(x) = 1, ~~ d F_t(x) = \langle x - a_t, A_t^{-1/2} d W_t \rangle F_t(x),
\end{equation}
$$
a_t = \int_{\RR^n} x f(x) F_t(x) dx, ~~ A_t = \int_{\RR^n} (x - a_t) \otimes (x - a_t) f(x) F_t(x) dx,
$$
for all $x \in \RR^n$ and all $0 \leq t \leq t_0$.
\end{lemma}
\emph{Proof:} \\
Fix $x \in \RR^n$. We will show that $d F_t(x) = \langle x - a_t, A_t^{-1/2} d W_t \rangle F_t(x)$. The correctness of the other equations is obvious.
Define, 
$$
G_t(x) = V_t F_t(x) = e^{\langle c_t, x \rangle - \frac{1}{2} \langle B_t x, x \rangle }.
$$
Equation (\ref{stochastic1}) clearly implies that $[B]_t = 0$. Let $Q_t(x)$ denote the quadratic variation of the process $\langle x, c_t \rangle$. We have,
$$
d \langle x, c_t \rangle = \langle A_t^{-1/2} x, dW_t + A_t^{-1/2} a_t dt \rangle.
$$
It follows that,
$$
d Q_t(x) = \langle A_t^{-1} x, x \rangle dt.
$$
Using It\^{o}'s formula, we calculate
$$
d G_t(x) = \left ( \langle x, d c_t \rangle - \frac{1}{2} \langle d B_t x, x \rangle + \frac{1}{2} d Q_t(x) \right ) G_t(x) =
$$
$$
\left ( \langle x, A_t^{-1/2} dW_t + A_t^{-1} a_t dt \rangle - \frac{1}{2} \langle A_t^{-1} x, x \rangle dt + \frac{1}{2} \langle A_t^{-1} x, x \rangle dt \right ) G_t(x) = 
$$
$$
\langle x, A_t^{-1/2} dW_t + A_t^{-1} a_t dt \rangle G_t(x).
$$
Next, we calculate,
$$
d V_t(x) = d \int_{\RR^n} e^{\langle c_t, x \rangle - \frac{1}{2} \langle B_t x, x \rangle }  f(x) dx = 
$$
$$
\int_{\RR^n} d G_t(x) f(x) dx = \int_{\RR^n} \langle x, A_t^{-1/2} dW_t + A_t^{-1} a_t dt \rangle G_t(x) f(x) dx = 
$$
$$
V_t \left \langle a_t, A_t^{-1/2} dW_t + A_t^{-1} a_t dt \right \rangle.
$$
So, using It\^{o}'s formula again,
$$
d V_t^{-1} = - \frac{d V_t}{V_t^2} + \frac{d [V]_t}{V_t^3} = 
$$
$$
- V_t^{-1} \left \langle a_t, A_t^{-1/2} dW_t + A_t^{-1} a_t dt \right \rangle + V_t^{-1} \langle A_t^{-1} a_t, a_t \rangle.
$$
Applying It\^{o}'s formula one last time yields,
$$
d F_t(x) = d (V_t^{-1} G_t(x)) = 
$$
$$
G_t(x) d V_t^{-1} + V_t^{-1} d G_t(x) + d [V^{-1}, G(x)]_t =
$$
$$
- V_t^{-1} \left \langle a_t, A_t^{-1/2} dW_t + A_t^{-1} a_t dt \right \rangle G_t(x) + V_t^{-1} \langle A_t^{-1} a_t, a_t \rangle G_t(x) + 
$$
$$
+ V_t^{-1} \langle x, A_t^{-1/2} dW_t + A_t^{-1} a_t dt \rangle G_t(x) - \langle A_t^{-1/2} a_t, A_t^{-1/2} x \rangle V_t^{-1} G_t(x) dt = 
$$
$$
\langle A_t^{-1/2} dW_t, x - a_t \rangle F_t(x).
$$
This finishes the proof.
\qed \\
\begin{remark}
In view of the above lemma it can be seen that, in some sense, the above is just the continuous version of the following iterative process: at every time step, multiply the function by a linear function equal to $1$ at the barycenter, whose gradient has a random direction distributed
uniformly on the ellipsoid of inertia. This construction may also be thought of as a variant of the Brownian motion on the Riemannian manifold constructed in \cite{EK1}.
\end{remark}

\begin{remark}
Rather than defining the process $F_t$ through equations (\ref{stochastic1}) and (\ref{defineF}), one may alternatively define it directly with the infinite system of stochastic differential equations in formula (\ref{contloc}). In this case, the existence and uniqueness of the solution can be shown using \cite[Theorem 5.2.2, page 159]{KX} (however, some extra work is needed in order to show that the conditions of this theorem hold).
\end{remark}
\medskip
In the remainder of this note, most of the calculations involving the process $f_t$ will use the formula (\ref{contloc}) rather
than the formulas (\ref{stochastic1}) and (\ref{defineF}). \\ \\
The remaining part of this section is dedicated to analyzing some basic properties of $\Gamma_t(f)$. We begin with:
\begin{lemma} \label{basic1}
The process $\Gamma_t(f)$ satisfies the following properties: \\ \\
(i) The function $\Gamma_t(f)$ is almost surely well defined, finite and log-concave for all $t > 0$. \\
(ii) For all $t > 0$, $\int_{\RR^n}f_t(x) dx = 1$. \\
(iii) The process has a semi-group property, namely, 
$$
\Gamma_{s+t}(f) \sim \frac{1}{\sqrt{\det A_s}} \Gamma_{t}(\sqrt{\det A_s} \Gamma_s(f) \circ L^{-1}) \circ L,$$
where
$$
L(x) = A_s^{-1/2}(x - a_s).
$$
(iv) For every $x \in \RR^n$, the process $f_t(x)$ is a martingale.
\end{lemma}
\bigskip
In order to prove (i), we will first need the following technical lemma:
\begin{lemma} \label{basic1.5}
For every dimension $n$, there exists a constant $c(n) > 0$ such that,
$$
\PP (A_{t} \geq c(n) Id, ~~ \forall 0 \leq t \leq c(n)) \geq c(n).
$$
\end{lemma}
The proof of this lemma is postponed section 6. \\ \\
\emph{Proof of lemma \ref{basic1}:} \\
To prove (i), we have to make sure that $A_t^{-1/2}$ does not blow up. To this end, define $t_0 = \inf \{t |~~ \det A_t = 0 \}$. By continuity, $t_0 > 0$. Equation (\ref{defineF}) suggests that $f_t$ is log-concave for all $t < t_0$. The fact that $t_0 = \infty$ will be proven below. \\
We start by showing that both (ii) and (iii) hold for any $t < t_0$. \\
We first calculate, using (\ref{contloc}),
$$
d \int_{\RR^n} f(x) F_t(x) dx = \int_{\RR^n} f(x) d F_t(x) dx =
$$
\begin{equation}
 \int_{\RR^n} f(x) F_t(x) \langle A_t^{-1/2} d W_t, x - a_t \rangle dx = 0,
\end{equation}
with probability 1. The last equality follows from the definition of $a_t$ as the barycenter of the measure $f(x) F_t(x) dx$. We conclude (ii). \\ \\
We continue with proving (iii). To do this, fix some $0 < s < t_0 - t$ and write,
\begin{equation} \label{normilization}
L(x) = A_s^{-1/2}(x - a_s).
\end{equation}
We normalize $f_s$ by defining,
$$
g(x) = \sqrt{\det A_s} f_s(L^{-1}(x)),
$$
which is clearly an isotropic probability density. Let us inspect $\Gamma_t(g(x))$. We have, using (\ref{contloc}),
$$
d \Gamma_t(g)(x) |_{t=0} = g(x) \langle x, d W_t \rangle =  \sqrt{\det A_s} f_s(L^{-1}(x)) \langle L (L^{-1}(x)), d W_t \rangle =
$$
$$
\sqrt{\det A_s} f_s(L^{-1}(x)) \langle L^{-1}(x) - a_s , A_s^{-1/2} d W_t \rangle,
$$
On the other hand,
$$
d f_s(L^{-1}(x)) = f_s(L^{-1}(x)) \langle L^{-1}(x) - a_s , A_s^{-1/2} d W_s \rangle
$$
in other words,
$$
d \Gamma_t(\sqrt{\det A_s} \Gamma_s(f) \circ L^{-1}) \left |_{t=0} \right . \sim \sqrt{\det A_s} d \Gamma_t(f) \circ L^{-1} \left |_{t=s} \right .
$$
which proves (iii). \\ \\
We are left with showing that $t_0 = \infty$.
To see this, write,
$$
s_1 = \min \{t ~~; ~~ ||A_t^{-1}||_{OP} = c^{-1}(n) \},
$$
where $c(n)$ is the constant from lemma \ref{basic1.5}. Note that, by continuity, $s_1$ is well-defined and almost-surely positive. When time $s$ comes, we may define $L_1$ as in (\ref{normilization}), and continue running the process on the function $f \circ L_1^{-1}$ as above. We repeat this every time $||A_t^{-1}||_{OP}$ hits the value $c^{-1}(n)$, thus generating the hitting times $s_1,s_2,...$.
Lemma \ref{basic1.5} suggests that,
$$
\PP \left . \left (s_{i+1} - s_i > c(n)~ \right | ~ s_1,s_2,...,s_i \right ) > c(n),
$$
which implies that, almost surely, $s_{i+1} - s_{i} > c(n)$ for infinitely many values of $i$. Thus, $\lim_{n \to \infty} s_n = \infty$ almost surely, and so $t_0 = + \infty$. \\
Part (iv) follows immediately from formula (\ref{contloc}). The lemma is proven. \qed \\ \\
Our next task is to analyze the path of the barycenter $a_t = \int_{\RR^n}x f_t(x) dx$. We have, using (\ref{contloc}),
\begin{equation} \label{pathbc}
d a_t = d \int_{\RR^n}x f(x) F_t (x) dx = \int_{\RR^n}x f(x) F_t(x) \langle x - a_t, A_t^{-1/2} d W_t \rangle dx =
\end{equation}
$$
\left ( \int_{\RR^n}(x - a_t) \otimes (x - a_t) f_t(x) dx \right ) (A_t^{-1/2} d W_t) = A_t^{1/2} d W_t.
$$
where the third equality follows from the defition of $a_t$, which implies,
$$
\int_{\RR^n}a_t f(x) F_t(x) \langle x - a_t, A_t^{-1/2} d W_t \rangle = 0.
$$
One of the crucial points, when using this localization scheme, will be to show that the barycenter of the measure does not move too much throughout the process. For this, we would like to attain upper bounds on the eigenvalues of the matrix $A_t$.
We start with a simple observation: \\ \\
Equation (\ref{defineF}) shows that the measure $f_t$ is log-concave with respect
to the measure $e^{- \frac{1}{2} |B_t^{1/2} x|^2}$. The following result, which is well-known to experts,
shows that measures which possess this property attain certain concentration inequalities.

\begin{proposition} \label{bakryemery}
There exists a universal constant $\Theta >0$ such that the following holds:
Let $\phi: \RR^n \to \RR$ be a convex function and let $K>0$.
Suppose that,
$$
d \mu(x) = Z e^{-\phi(x) - \frac{1}{2 K^2} |x|^2} dx
$$
is a probability measure whose barycenter lies at the origin. Then, \\
(i) For all Borel sets
$A \subset \RR^n$, with $0.1 \leq \mu(A) \leq 0.9$, one has,
$$
\mu(A_{K \Theta}) \geq 0.95
$$
where $A_{K \Theta}$ is the $K \Theta$-extension of $A$, defined in the previous section. \\
(ii) For all $\theta \in \Sph$,
$$
\int \langle x, \theta \rangle^2 d \mu(x) \leq \Theta K^2.
$$
\end{proposition}
\emph{Proof:} \\
Denote the density of $\mu$ by $\rho(x)$. Let $B$ be the complement of $A_{K \Theta}$, where the constant $\Theta$
will be chosen later on. Define,
$$
f(x) = \rho(x) \mathbf{1}_{A}, ~~ g(x) = \rho(x) \mathbf{1}_{B}.
$$
Note that for $x \in A$ and $y \in B$, we have $|x-y|>K \Theta$. Thus, by the parallelogram law,
$$
\left | \frac{x+y}{2} \right |^2 \leq \frac{|x|^2 + |y|^2}{2} - \frac{1}{4} K^2 \Theta^2,
$$
which implies,
$$
e^{-\frac{1}{2 K^2} \left |\frac{x+y}{2} \right |^2} \geq \sqrt{e^{-\frac{1}{2K^2} |x|^2} e^{-\frac{1}{2K^2} |x|^2} } e^{\frac{1}{8} \Theta^2}.
$$
Since the function $\phi$ is assumed to be convex, we obtain
$$
\rho \left( \frac{x+y}{2} \right ) \geq \sqrt{f(x) g(y)}e^{\frac{1}{8} \Theta^2}.
$$
Now, using the Prekopa-Leindler theorem, we attain
$$
\mu(A) \mu(B) = \int_{\RR^n} f(x) dx \int_{\RR^n} g(x) dx \leq e^{-\frac{1}{4} \Theta^2},
$$
so,
$$
\mu(A_{K \Theta}) \geq 1 - \frac{e^{-\frac{1}{4} \Theta^2}}{\mu(A)} \geq 1 - 10 e^{-\frac{1}{4} \Theta^2}.
$$
Clearly, a large enough choice of the constant $\Theta$ gives (i). To prove (ii), we define,
$$
g(t) = \mu(\{x; \langle x, \theta \rangle \geq t \}),
$$
and take $A = \{x; \langle x, \theta \rangle < g^{-1}(0.5) \}$. An application on (i) on the set $A$ gives,
$$
g(g^{-1}(0.5) + K \Theta) \leq 0.05
$$
since $g$ is log-concave, we attain
$$
g(g^{-1}(0.5) + t K \Theta) \leq 0.05^t, ~~ \forall t > 1
$$
and in the same way, one can attain,
$$
1 - g(g^{-1}(0.5) - t K \Theta) \leq 0.05^t, ~~ \forall t > 1.
$$
Part (ii) of the proposition is a direct consequence of the last two equations. \qed \\ \\
Plugging (\ref{defineF}) into part (ii) of this theorem gives,
\begin{equation} \label{goodbound}
A_t \leq \Theta ||B_t^{-1}||_{OP} Id \leq \Theta \left( \int_0^t  \frac{ds}{||A_s||_{OP}} \right )^{-1} Id, ~~~ \forall t>0.
\end{equation}
By our assumption (\ref{compactsupport}) we deduce that $A_t$ is bounded by $n^2 Id$, which immediately gives
\begin{equation} \label{poorbound}
A_t < \frac{\Theta n^2}{t} Id.
\end{equation}
$$
~
$$
The bound (\ref{poorbound}) will be far from sufficient for our needs, and the next section is dedicated
to attaining a better upper bound. However, it is good enough to show that the barycenter, $a_t$, converges in distribution to the density $f(x)$. \\ \\
Indeed, (\ref{poorbound}) implies that 
\begin{equation} \label{atconverges}
\lim_{t \to \infty} W_2(f_t, \delta_{a_t}) = 0
\end{equation}
where $\delta_{a_t}$ is the probability measure supported on $\{a_t\}$. In other words the probability density $f_t(x)$ converges to a delta measure. By the martingale property, part (iv) of lemma \ref{basic1}, we know that $\EE[f_t(x)] = f(x)$, thus, $X_t := a_t$ converges, in Wasserstein metric, to the original random vector $X$ as $t \to \infty$.
\begin{remark}
It is interesting to compare this construction with the construction by Lehec in \cite{Lehec}. In both cases, a certain It\^{o} process converges
to a given log-concave measure. In the result of Lehec, the convergence is ensured by applying a certain adapted \emph{drift}, while
here, it is ensured by adjusting the \emph{covariance matrix} of the process.
\end{remark}

We end this section with a simple calculation in which we analyze the process $\Gamma_t(f)$ in the simple case that $f$ is the standard Gaussian measure. While the calculation will not be necessary for our proofs, it may provide the reader a better understanding of the process. Define,
$$
f(x) = (2 \pi)^{-n/2} e^{-|x|^2 / 2}. 
$$
According to formula (\ref{defineF}), the function $f_t$ takes the form,
$$
f_t(x) = V_t^{-1} \exp \left ( \langle x, c_t \rangle - \frac{1}{2} \left \langle (B_t + Id) x, x \right \rangle \right )
$$
where $V_t \in R, c_t \in \RR^n$ are certain It\^{o} processes. It follows that the covariance matrix $A_t$ satisfies,
$$
A_t^{-1} = B_t + Id.
$$
Recall that $B_t = \int_0^t A_s^{-1} ds$. It follows that,
$$
\frac{d}{dt} B_t = B_t + Id, ~~ B_0 = 0.
$$
So,
$$
B_t = (e^t - 1) Id,
$$
which gives,
$$
A_t = e^{-t} Id.
$$
Next, we use (\ref{pathbc}) to derive that,
$$
d a_t = e^{-t/2} d W_t,
$$
which implies,
$$
a_t \sim W_{1 - \exp(-t)}.
$$
We finally get,
$$
f_t = e^{nt/2} (2 \pi)^{-n/2} \exp \left (-\frac{1}{2} e^t \left |(x - W_{1 - \exp(-t)}) \right |^2 \right ).
$$
\section{Analysis of the matrix $A_t$}

In the previous section we saw that the covariance matrix of the $f_t$, $A_t$, satisfies (\ref{poorbound}).
The goal of this section is to give a better bound, which holds also for small $t$. Namely, we want to prove:
\begin{proposition} \label{mainsec3}
There exist universal constants $C,c > 0$ such that for all $n \geq 2$ the following holds: Let $f:\RR^n \to \RR_+$ be an isotropic, log concave probability density. Let $A_t$ be the covariance matrix of $\Gamma_t(f)$. Then, \\
(i) Define the event $F$ by,
\begin{equation} \label{deff}
F := \left \{ ||A_t||_{OP} < C K_n^2 (\log n) e^{-c t}, ~~ \forall t > 0 \right \}.
\end{equation}
One has,
\begin{equation} \label{onb1}
\PP(F) ~~\geq~~ 1 - (n^{-10}).
\end{equation}
(ii) For all $t > 0$, $\EE[Tr(A_t)] \leq n$. \\
(iii) Whenever the event $F$ holds, the following also holds: \\
For all $t > \frac{1}{K_n^2 \log n}$ there exists a convex function $\phi_t(x)$ such that the function $f_t$ is of the form,
\begin{equation} \label{ftForm2}
f_t(x) = \exp \left ( - \left | \frac{x}{C K_n \sqrt{\log n}} \right |^2 - \phi_t(x) \right).
\end{equation}
\end{proposition}
\bigskip
Before we move on to the proof, we will establish some simple properties of the matrix $A_t$. Our first task is
to find the differential of process $A_t$. We have, using It\^{o}'s formula with equation (\ref{contloc}),
\begin{equation} \label{dAt1}
d A_t = d \int_{\RR^n}(x - a_t) \otimes (x - a_t) f_t(x) dx =
\end{equation}
$$
\int_{\RR^n}(x - a_t) \otimes (x - a_t) d f_t(x) dx - 2 \int_{\RR^n} d a_t \otimes (x - a_t) f_t(x) dx -
$$
$$
- 2 \int_{\RR^n}(x - a_t) \otimes d [a_t, f_t(x)]_t dx + d [a_t, a_t] \int_{\RR^n} f_t(x) dx.
$$
Let us try to understand each of this terms. The second term is,
$$
\int_{\RR^n} d a_t \otimes (x - a_t) f_t(x) dx = d a_t \otimes \int_{\RR^n} (x - a_t) f_t dt = 0.
$$
Recall that by (\ref{pathbc}), $d a_t = A_t^{1/2} dW_t$, which gives,
\begin{equation} \label{dAt2}
d [a_t, a_t]_t = A_t dt,
\end{equation}
and
$$
d[a_t, f_t(x)] = f_t(x) A_t^{1/2} A_t^{-1/2} x dt = f_t(x) x dt,
$$
which implies,
\begin{equation} \label{dAt3}
\int_{\RR^n} (x - a_t) \otimes d [a_t, f_t(x)]_t dx = \int_{\RR^n} (x - a_t) \otimes x f_t(x) dx dt = 
\end{equation}
$$
\int_{\RR^n} (x - a_t) \otimes (x - a_t) f_t(x) dx dt = A_t dt
$$
Plugging equations (\ref{dAt1}), (\ref{dAt2}) and (\ref{dAt3}) together gives,
$$
d A_t = \int_{\RR^n}(x - a_t) \otimes (x - a_t) d f_t(x) dx - A_t dt,
$$
so we finally get,
$$
d A_t = \int_{\RR^n}(x - a_t) \otimes (x - a_t) \langle x - a_t, A_t^{-1/2} dW_t \rangle f_t(x) dx - A_t dt.
$$
Note that the term $A_t dt$ is positive definite, hence, the fact it appears in the differential can only make all of the eigenvalues of $A_t$ smaller (as a matter of fact, this term induces a rather strong drift of all the eigenvalues towards $0$, which we will not even use). Consequently, we can define $\tilde A_t = A_t + \int_0^t A_s ds$, so that
\begin{equation} \label{datilde}
d \tilde A_t = \int_{\RR^n}(x - a_t) \otimes (x - a_t) \langle x - a_t, A_t^{-1/2} d W_t \rangle f_t(x) dx
\end{equation}
and $\tilde A_0 = A_0 = Id$. Clearly, $A_t \leq \tilde A_t$ for all $t>0$. In order to control $||A_t||_{OP}$, it is thus enough to bound $||\tilde A_t||_{OP}$. \\ \\
For a fixed value of $t$, let $v_1,...,v_n$ be an orthonormal basis, with respect to which $\tilde A_t$ is diagonal, and write
$\alpha_{i,j} = \langle v_i, \tilde A_t v_j \rangle $ for the entries of $\tilde A_t$ with respect to this basis. Equation (\ref{datilde}) can be written,
$$
d \alpha_{i,j} = \int_{\RR^n} \langle x, v_i \rangle \langle x, v_j \rangle \langle A_t^{-1/2} x, d W_t \rangle f_t (x + a_t) dx.
$$
Next, denote
\begin{equation} \label{defxi}
\xi_{i,j} = \frac{1}{\sqrt{\alpha_{i,i} \alpha_{j,j}}} \int_{\RR^n} \langle x, v_i \rangle \langle  x, v_j \rangle A_t^{-1/2} x f_t (x + a_t) dx.
\end{equation}
So,
\begin{equation} \label{dalpha}
d \alpha_{i,j} = \sqrt{\alpha_{i,i} \alpha_{j,j}} \langle \xi_{i,j}, d W_t \rangle, 
\end{equation}
and
$$
\frac{d}{dt} [\alpha_{i,j}]_t = \alpha_{i,i} \alpha_{j,j} |\xi_{i,j}|^2.
$$
$$
~
$$
As we will witness later, behaviour norm of the matrix $\tilde A_t$ depends highly on the norms of the vectors $\xi_{i,j}$,
which induce a certain repulsion between the eigenvalues. The next lemma will come in handy when we need to bound these norms:

\begin{lemma} \label{xibounds}
The vectors $\xi_{i,j}$ satisfy the following bounds: \\
(i) For all $1 \leq i \leq n$, $|\xi_{i,i}| < C$ for some universal constant $C>0$. \\
(ii) For all $1 \leq i \leq n$, $\sum_{j=1}^n | \xi_{i,j} |^2 \leq K_n^2$.
\end{lemma}
\emph{Proof:} \\
Since $A_t^{1/2} v_i = \sqrt{\alpha_{i,i}} v_i$ for all $1 \leq i \leq n$, we have
$$
\xi_{i,j} = \int_{\RR^n} \langle A_t^{-1/2} x, v_i \rangle \langle  A_t^{-1/2} x, v_j \rangle A_t^{-1/2} x f_t (x + a_t) dx
$$
Define again, as above, $\tilde f_t(x) = \sqrt {\det A_t} f_t(A_t^{1/2} x + a_t)$. 
By substituting $y = A_t^{-1/2} x$, the equation becomes
\begin{equation} \label{defxigood}
\xi_{i,j} = \int_{\RR^n} \langle y, v_i \rangle \langle  y, v_j \rangle y \tilde f_t (y) dy.
\end{equation}
Recall that $\tilde f_t$ is isotropic. The last equation shows that the vectors $\xi_{i,j}$, in some sense, do not depend on the
position of $f$. Using the Cauchy-Schwartz inequality, one has
\begin{equation} \label{estxi1}
|\xi_{i,i}| = \left | \int_{\RR^n} \langle y, v_i \rangle^2 y \tilde f_t(y) dy \right | = \left | \int_{\RR^n} \langle y, v_i \rangle^2 \left \langle y, \frac{\xi_{i,i}} {|\xi_{i,i}|} \right \rangle \tilde f_t(y) dy \right | \leq
\end{equation}
$$
\sqrt{\int_{\RR^n} \langle y, v_i \rangle^4 \tilde f_t(y) dy
\int_{\RR^n}  \left \langle y, \frac{\xi_{i,i}} {|\xi_{i,i}|} \right \rangle^2 \tilde f_t(y) dy}.
$$
A well-known fact about isotropic log-concave measures (see for example \cite[Lemma 5.7]{LV}) is that for every
$p>0$ there exists a constant $c(p)$ such that for every isotropic log-concave density $\rho(x)$ on $\RR^n$ and every $\theta \in \Sph$,
$$
\int_{\RR^n} |\langle x, \theta \rangle|^p \rho(x) dx \leq c(p).
$$
Using this with (\ref{estxi1}) establishes (i). Next, by the definition of $K_n$, we have for all $1 \leq i \leq n$,
$$
\sum_{j=1}^n |\xi_{i,j}|^2 = \sum_{j=1}^n \sum_{k=1}^n \left | \int_{\RR^n} \langle y, v_i \rangle \langle  y, v_j \rangle \langle y, v_k \rangle \tilde f_t (y) dy \right |^2 = 
$$
$$
\left | \left | \int_{\RR^n} y \otimes y \langle y, v_i \rangle \tilde f_t (y) dy  \right| \right |_{HS}^2 \leq K_n^2.
$$
The lemma is proven. \qed \\ \\
We are now ready to prove the main proposition of the section. \\ \\
\emph{Proof of proposition \ref{mainsec3}}: \\
We fix a positive integer $p$ whose value will be chosen later, and define,
\begin{equation} \label{defst}
S_t = Tr \left (\tilde A_t^p \right).
\end{equation}
Since $S_t$ is a smooth function of the coefficients $\{ \alpha_{i,j} \}$, which are It\^{o} processes (assuming that the basis $v_1,...,v_n$ is fixed), $S_t$ itself is also an It\^{o} process. Fix some $t>0$. Our next goal will be to find $d S_t$. To that end, define $\Gamma$ to be the set of $(p+1)$-tuples, $(j_1,,..,j_{p+1})$, such that $j_i \in \{1,...,n\}$ for all $1 \leq i \leq p+1$ and such that $j_1=j_{p+1}$. It is easy to verify that,
\begin{equation} \label{stpaths}
S_t = \sum_{(j_1,...,j_{p+1}) \in \Gamma} \alpha_{j_1,j_2} \alpha_{j_2,j_3} \cdots \alpha_{j_p,j_{p+1}}.
\end{equation}
Since $Tr(\tilde A_t^p)$ does not depend on the choice of orthogonal coordinates, after fixing the value of $t$, we are free to choose our coordinates such that the matrix $\tilde A_t$ is diagonal, thus assuming that $\alpha_{i,j}=0$ whenever $i \neq j$ and that
(\ref{dalpha}) holds (in other words, we calculate the differential $d S_t$ using a basis $v_1,...,v_n$ which depends on $t$. However, after fixing the value of $t$, the calculation itself is with respect to a fixed basis). A moment of reflection reveals that, in this case, the term $d (\alpha_{j_1,j_2} \cdots \alpha_{j_p,j_{p+1}})$ can be non-zero only if there are at most two distinct indices $i_1, i_2$ such that $j_{i_1} \neq j_{i_1 + 1}$ and $j_{i_2} \neq j_{i_2 + 1}$. We are left with two types of terms whose differential is non-zero. The first type of term contains no off-diagonal entries, and has the form $(\alpha_{i,i})^p$. 
Using equation (\ref{dalpha}), we calculate its differential,
\begin{equation} \label{term1est}
d \alpha_{i,i}^p = p \alpha_{i,i}^{p-1} d \alpha_{i,i} + \frac{1}{2} p(p-1) \alpha_{i,i}^{p-2} d [\alpha_{i,j}]_t =
\end{equation}
$$
p \alpha_{i,i}^{p}  \langle \xi_{i,i}, d W_t \rangle + p(p-1) \alpha_{i,i}^p |\xi_{i,i}|^2 dt.
$$
The second type of term will contain exactly two off-diagonal entries, and due to the symmetry of the matrix and
the constraint $j_{1} = j_{p+1}$, it has the form: 
$$
(\alpha_{i,i})^{k_1} \alpha_{i,j} (\alpha_{j,j})^{k_2} \alpha_{j,i} (\alpha_{i,i})^{k_3} =
(\alpha_{i,i})^k  (\alpha_{j,j})^{p-k-2} (\alpha_{i,j})^{2}
$$
where $i \neq j$ and $0 \leq k \leq p-2$. Keeping in mind that $\alpha_{i,j} = 0$, we calculate,
$$
d \left ( (\alpha_{i,i})^k  (\alpha_{j,j})^{p-k-2} (\alpha_{i,j})^{2}  \right ) = (\alpha_{i,i})^k  (\alpha_{j,j})^{p-k-2} \left (2 \alpha_{i,j} d \alpha_{i,j} + d [\alpha_{i,j}]_t \right ) =
$$
$$
(\alpha_{i,i})^{k+1} (\alpha_{j,j})^{p-k-1} |\xi_{i,j}|^2 dt.
$$
We may clearly assume $\alpha_{1,1} \geq \alpha_{2,2} \geq ... \geq \alpha_{n,n}$, which implies that for $i < j$ and for all values
of $k$, one has
\begin{equation} \label{term2est}
d \left ( (\alpha_{i,i})^k  (\alpha_{j,j})^{p-k-2} (\alpha_{i,j})^{2} \right ) \leq (\alpha_{i,i})^p |\xi_{i,j}|^2 dt
\end{equation}
Inspect the equation (\ref{stpaths}). For every $1 \leq i \leq n$, the expansion on the right hand side contains
exactly one term of the first type, and for every distinct $i,j$ with $i \neq j$, it contains $\frac{p(p-1)}{2}$ terms
of the second type (or otherwise, for all choices such that $i < j$, it contains $p(p-1)$ terms of this type). Using (\ref{term1est}) and (\ref{term2est}), we conclude
$$
d S_t \leq \sum_{i=1}^n p \alpha_{i,i}^{p}  \langle \xi_{i,i}, d W_t \rangle + p(p-1) \alpha_{i,i}^p |\xi_{i,i}|^2 dt +
\sum_{1 \leq i, j \leq n \atop i < j} p(p-1) (\alpha_{i,i})^p |\xi_{i,j}|^2 dt \leq 
$$
$$
\sum_{i=1}^n p \alpha_{i,i}^{p} \langle \xi_{i,i}, d W_t \rangle + p^2 \sum_{i=1}^n \alpha_{i,i}^p \sum_{j=1}^n |\xi_{i,j}|^2 dt \leq
$$
$$
\sum_{i=1}^n p \alpha_{i,i}^{p} \langle \xi_{i,i}, d W_t \rangle + p^2 S_t K_n^2 dt
$$
where in the last inequality we used the part (ii) of lemma \ref{xibounds}. \\ \\
A well-known property of It\^o processes is existence and uniqueness of the decomposition $S_t = M_t + E_t$, where $M_t$ is a local martingale and $E_t$ is an adapted process of locally bounded variation. In the last equation, we attained,
\begin{equation} \label{dst}
d E_t \leq p^2 K_n^2 S_t dt,
\end{equation}
and also,
$$
\frac{d [S]_t}{dt} = \left | \sum_{i=1}^n p \alpha_{i,i}^{p} \xi_{i,i} \right |^2.
$$
Using part (i) of lemma \ref{xibounds} yields,
\begin{equation} \label{dst2}
\frac{d [S]_t}{dt} \leq C p^2 S_t^2.
\end{equation}
Next, we use the unique decomposition $\log S_t = Y_t + Z_t$ where $Y_t$ is a local martingale, $Z_t$ is an adapted process of locally bounded variation and $Y_0 = 0$. According to It\^{o}'s formula and formula (\ref{dst2}),
\begin{equation} \label{dyt}
\frac{d [Y]_t}{dt} = \frac{1}{S_t^2} \frac{d [S]_t}{dt} \leq C p^2.
\end{equation}
By Dambis / Dubins-Schwartz theorem, we know that there exists a standard Wiener process
$\tilde W_t$ such that $Y_t$ has the same distribution as $\tilde W_{[Y]_t}$. An application of the so-called
reflection principle gives,
$$
\PP \left (\max_{t \in [0, p] } \tilde W_t \geq t p \right ) = 
$$
$$
2 \PP (\tilde W_p \geq t p ) < C e^{-\frac{1}{2} t^2 p }.
$$
Choosing $t$ to be a large enough universal constant, $C_1$, yields
$$
\PP \left (\max_{t \in [0, p] } \tilde W_t \geq C_1 p \right ) < e^{-10 p},
$$
(where we used the fact that $p \geq 1$). Using (\ref{dyt}), we attain
$$
\PP \left (\max_{t \in \left [0, \frac{1}{p} \right ] } Y_t > C_2 p \right ) < e^{- 10 p}
$$
for some universal constant $C_2 > 0$. We now use It\^{o}'s formula again, this time with formula (\ref{dst}), to get
$$
\frac{d}{dt} Z_t = \frac{1}{S_t} \frac{d}{dt} E_t - \frac{1}{2 S_t^2} \frac{d [S]_t}{dt} \leq K_n^2 p^2.
$$
The last two equations and the legitimate assumption that $K_n \geq 1$ give,
$$
\PP \left (\max_{t \in \left [0, \frac{1}{K_n^2 p} \right ] } \log S_t - \log n > C p \right ) < e^{- 10 p}.
$$
We choose $p = \lceil \log n \rceil$ to get,
$$
\PP \left (\max_{t \in \left [0, \frac{1}{K_n^2 \log n} \right ]} S_t^{1 / \lceil \log n \rceil} > C' \right ) < \frac{1}{n^{10}},
$$
for some universal constant $C'>0$.
Define the event $F$ as the complement of the event in the equation above,
$$
F := \left \{\max_{t \in \left [0, \frac{1}{K_n^2 \log n} \right ]} S_t^{1/\lceil \log n \rceil} \leq C' \right \}.
$$
Clearly, whenever the event $F$ holds, we have,
\begin{equation} \label{wheneholds}
||A_t||_{OP} \leq ||\tilde A_t||_{OP} \leq C', ~~~ \forall t \in \left [0, \frac{1}{K_n^2 \log n} \right ].
\end{equation}
Our next task is to bound the norm for larger values of $t$. To this end, recall the bound (\ref{goodbound}). Recalling that $B_t = \int_0^t
A_s^{-1} ds$, and applying (\ref{goodbound}) gives,
$$
\frac{d}{dt} B_t = A_t^{-1} \geq \frac{Id}{\Theta ||B_t^{-1}||_{OP}}.
$$
So,
\begin{equation} \label{btexp1}
\frac{d}{dt} \frac{1}{||B_t^{-1}||_{OP}} \geq \frac{1}{\Theta ||B_t^{-1}||_{OP}}.
\end{equation}
By the definition of $B_t$ and by (\ref{wheneholds}), it follows that whenever $F$ holds one has,
\begin{equation} \label{btexp2}
\frac{1}{||B_{\delta^2}^{-1}||_{OP}} \geq C \delta^2
\end{equation}
where $\delta^2 = \frac{1}{K_n^2 \log n}$.
Equations (\ref{btexp1}) and (\ref{btexp2}) imply,
$$
B_{t} \geq c \delta^2 e^{(t - \delta^2) / \Theta } Id, ~~~ \forall t > \delta^2
$$
which gives, using (\ref{goodbound}),
$$
A_{t} \leq C \delta^{-2} e^{(\delta^2 - t) / \Theta } Id.
$$
Part(i) of the proposition is established. In order to prove the bound for $\EE[Tr(A_t)]$, write
$S_t = \sum_{i=1}^n Tr(\tilde A_t)$. Setting $p = 1$ in (\ref{defst}) gives,
$\frac{d}{dt} \EE[S_t] = 0$, which implies (ii). Part (iii) of the proposition follows directly from equations (\ref{btexp2}) and (\ref{defineF}). The proposition is complete.
\qed \\

Proposition \ref{bakryemery} gives an
immediate corollary to part (iii) of proposition \ref{mainsec3}:
\begin{corollary} \label{goodcheeger}
There exist universal constants $c, \Theta >0$ such that whenever the event $F$ defined in (\ref{deff}) holds, the
following also holds: \\
Define $\delta = \frac{1}{K_n \sqrt{\log n}}$. Let $t > \delta^2$ and let $E \subset \RR^n$ be a measurable set
which satisfies,
\begin{equation} \label{defR}
0.1 \leq \int_E f_t(x) dx \leq 0.9.
\end{equation}
One has,
\begin{equation}
\int_{E_{\Theta / \delta } \setminus E} f_t(x) dx \geq c
\end{equation}
where $E_{\Theta / \delta}$ is the $\frac \Theta \delta$-extension of $E$, defined in the introduction.
\end{corollary}

\section{Thin shell implies spectral gap}

In this section we use the localization scheme constructed in the previous sections in order to prove
theorem \ref{mainthm1}.\\ \\
Let $f(x)$ be an isotropic log-concave probability density in $\RR^n$ and let $E \subset \RR^n$ be a measurable set.
Suppose that,
\begin{equation}
\int_E f(x) dx = \frac{1}{2}.
\end{equation}
Our goal in this section is to show that,
\begin{equation} \label{needtoshow11}
\int_{E_{\Theta / \delta} \setminus E} f(x) dx \geq c
\end{equation}
for some universal constants $c, \Theta >0$, where $\delta = \frac{1}{K_n \sqrt{\log n}}$ and $E_{\Theta / \delta}$ is the $\frac \Theta \delta$-extension of $E$. \\ \\
The idea is quite simple. Define $f_t := \Gamma_t(f)$, the localization of $f$ constructed in section 2, and fix $t > 0$. By the martingale property of the localization, we have,
\begin{equation} \label{martlocE}
\int_{E_{\Theta / \delta} \setminus E} f(x) dx = \EE \left [ \int_{E_{\Theta / \delta} \setminus E} f_t(x) dx  \right ].
\end{equation}
Corollary \ref{goodcheeger} suggests that if $t$ is large enough, the right term can be bounded from below if we only manage to bound the integral $\int_E f_t(x) dx$ away from 0 and from 1. \\ \\
Define,
$$
g(t) = \int_E f_t(x) dx.
$$
\bigskip
In view of the above, we would like to prove:
\begin{lemma} \label{gtbound}
There exists a universal constant $T>0$ such that,
$$
\PP \left (0.1 \leq g \left ( t \right ) \leq 0.9 \right ) > 0.5, ~~~ \forall t \in [0,T].
$$
\end{lemma}
\emph{Proof:} \\
We calculate, using (\ref{contloc}),
\begin{equation} \label{dgt}
d g(t) = \int_E f_t(x) \langle x - a_t, A_t^{-1/2} dW_t \rangle dx =
\end{equation}
(substitute $y = A_t^{-1/2}(x - a_t)$)
$$
\sqrt{\det A_t} \int_{A_t^{-1/2}(E - a_t)} f_t(A_t^{1/2} y + a_t) \langle y, d W_t \rangle dy =
$$
$$
\left \langle \sqrt{\det A_t} \int_{A_t^{-1/2}(E - a_t)} f_t(A_t^{1/2} y + a_t) y dy , d W_t \right \rangle.
$$
Define,
$$
\tilde f_t = \sqrt{\det A_t} f_t(A_t^{1/2} y + a_t), ~~~E_t = A_t^{-1/2} (E - a_t)
$$
The above equation becomes,
\begin{equation} \label{dgt2}
d g(t) = \left \langle \int_{E_t} y \tilde f_t(y) dy, d W_t \right \rangle.
\end{equation}
Assume, for now, that $\int_{E_t} y \tilde f_t(y) dy \neq 0$ and define $\theta = \frac{\int_{E_t} y \tilde f_t(y) dy}{|\int_{E_t} y \tilde f_t(y) dy|}$. Observe that, by definition, $\tilde f_t$ is isotropic. Consequently,
$$
\left |\int_{E_t} y \tilde f_t(y) dy \right| = \left | \int_{E_t} \langle y, \theta \rangle \tilde f_t(y) dy  \right | \leq
$$
$$
\int_{E_t} |\langle y, \theta \rangle| \tilde f_t(y) dy \leq \sqrt{ \int_{E_t} \langle y, \theta \rangle^2 \tilde f_t(y) dy} \leq 1.
$$
We therefore learn that,
$$
\frac{d}{dt} [g]_t \leq 1, ~~ \forall t>0.
$$
Define $h(t) = (g(t) - 0.5)^2$. By It\^{o}'s formula,
$$
d h(t)  = 2 (g(t) - 0.5) d g(t) + d [g]_t.
$$
Plugging the last two equations together gives,
$$
E[(g(t) - 0.5)^2] \leq t.
$$
The lemma follows from an application of Chebyshev's inequality. \qed \\ \\
The last ingredient needed for our proof is a theorem of E. Milman, \cite[Theorem 2.1]{Mil2}. The following
is a weaker formulation of this theorem which will be suitable for us: \\
\begin{theorem}[E. Milman] \label{thmmilman}
Suppose that a log-concave probability measure $\mu$ satisfies the following: there exist two constants, $0 < \lambda < \frac{1}{2}$ and
$\Theta > 0$, such that for all measurable $E \subset \RR^n$ with $\mu(E) \geq \frac{1}{2}$, one has 
$\mu(E_{\Theta}) \geq 1 - \lambda$. In this case, the measure $\mu$ satisfies the following isoperimetric inequality: \\
For all measurable $E \subset \RR^n$ with $\mu(E) \leq \frac{1}{2}$,
\begin{equation} \label{milman}
\frac{\mu^+(E)}{\mu(E)} \geq \frac{1 - 2 \lambda}{\Theta}.
\end{equation}
\end{theorem}
\medskip
Note that equation (\ref{milman}) is the exact type of inequality defining the constant $G_n$ in equation (\ref{defg}). We are now ready to prove the main proposition of this section. \\ \\
\emph{Proof of proposition \ref{gk}}: \\
Let $T$ be the constant from lemma \ref{gtbound}. Denote,
$$
G = \left \{ 0.1 \leq g(T) \leq 0.9 \right \} \cap F.
$$
where $F$ is the event defined in (\ref{deff}). According to lemma \ref{gtbound} and to (\ref{onb1}), one has $\PP(G) > 0.4$ for all $n \geq 2$. \\ \\
By (\ref{martlocE}) and by corollary \ref{goodcheeger}, there exist universal constants $\tilde c, \Theta > 0$ such that
\begin{equation}
\int_{E_{\Theta / \delta} \setminus E} f(x) dx = \EE \left [ \int_{E_{\Theta / \delta} \setminus E} f_T(x) dx  \right ] \geq
\end{equation}
$$
P(G) \EE \left [ \left . \int_{E_{\Theta / \delta} \setminus E} f_T(x) dx  ~~ \right | G \right ] \geq \tilde c.
$$
The result now follows directly from an application of theorem \ref{thmmilman}.
\qed \\
\begin{remark}
In the above proof, we used E. Milman's result in order to reduce the theorem to the case where $\int_E f(x) dx$ is exactly $\frac{1}{2}$, as well as to attain an isoperimetric inequality from a certain concentration inequality for distance functions. Alternatively, we may have replaced propsition \ref{bakryemery} with an essentially stronger result due to Bakry-Emery, proven in \cite{BE} (see also Gross, \cite{Gross}). Their result, which relies on the \emph{hypercontractivity principle}, asserts that a density of the form (\ref{defineF}) actually possesses a respective Cheeger constant. Using this fact, we may have directly bounded from below the surface area of any set with respect to the measure whose density is $f_t$.
\end{remark}
\medskip
The proof of lemma \ref{ksigma} is in section 6. Along with this lemma, we have established theorem
\ref{mainthm1}.

\section{ Stability of the Brunn-Minkowski inequality }
The main goal of this section is to prove theorem \ref{mainthm2}. \\ \\
The idea of the proof is as follows: Given two log-concave densities, $f$ and $g$,
we run the localization process we constructed in section 2 on both functions, so that their corresponding localization processes
are coupled together in the sense that we take the same Wiener process $W_t$ for both functions.
Recall formula (\ref{atconverges}), whose point is that the barycenters of the localized functions $f_t$ and $g_t$ converge, in the Wasserstein metric,
to the measures whose densities are $f$ and $g$, respectively. In view of this, it is enough to consider the paths of the barycenters and show that
they remain close to each other along the process. Recall that if $a_t$ is the barycenter of $f_t$, we have
$d a_t = A_t^{1/2} d W_t$. This formula tells us that as long as we manage to keep the covariance matrices of
$f_t$ and $g_t$  approximately similar to each other, the barycenters will not move too far apart. In order to do this, we use an idea from \cite{EK2}: when the integral of the supremum convolution of two given densities is rather small, these densities can essentially be regarded as parallel sections of an isotropic convex body, which means, by thin-shell concentration, that the corresponding covariance matrices cannot be very different from each other. \\ \\
We begin with some notation. For two functions $f,g: \RR^n \to \RR^+$, denote by $H_\lambda (f,g)$ the supremum convolution of the two functions, hence,
$$
H (f,g)(x) := \sup_{y \in \RR^n} \sqrt{f (x + y) g (x - y)}.
$$
Next, define,
$$
K(f,g) = \int_{\RR^n} H (f,g)(x)  dx.
$$
$$
~
$$
The following lemma is a variant of lemma 6.5 from \cite{EK2}.

\begin{lemma} \label{transportation}
There exists a universal constant $C>0$ such that the following holds: Let $f,g$ be log-concave probability densities in $\RR^n$. Define,
$$
A = Cov(f)^{-1/2} Cov(g) Cov(f)^{-1/2} - Id,
$$
and let $\{\delta_i \}_{i=1}^n$ be the eigenvalues of $A$ such that the order of 
$|\delta_i - 1|$ is decreasing. Then,
\begin{equation} \label{smallis}
|\delta_i - 1| \leq C K(f,g)^4, ~~ \forall 1 \leq i \leq n
\end{equation}
and,
\begin{equation} \label{largeis}
|\delta_i - 1| \leq C K(f,g) \tau_n i^{\kappa - \frac{1}{2}}, ~~\forall (\log K(f,g))^{C_1} \leq i \leq n
\end{equation}
where $C,C_1 > 0$ are universal constants.
\end{lemma}
\bigskip
Our main ideas in this section are contained in the following lemma: \\
\begin{lemma} \label{mainlemmasec5}
Let $\epsilon > 0$ and let $f$, $g$ be log-concave probability densities in $\RR^n$ such that $f$ is isotropic and the barycenter of $g$ lies at the origin. In that case, there exist two densities, $\tilde f, \tilde g$, which satisfy,
$$
\tilde f(x) \leq f(x), ~~\tilde g(x) \leq g(x), ~~ \forall x \in \RR^n,
$$
$$
\int_{\RR^n} \tilde f(x) dx = \int_{\RR^n} \tilde g(x) dx \geq 1 - \epsilon
$$
and,
\begin{equation} \label{smallw2}
W_2(\tilde f, \tilde g) \leq \frac{C}{\epsilon^{6}} \tau_n K(f,g)^{5} n^{2 (\kappa - \kappa^2) + \epsilon}
\end{equation}
\end{lemma}

\emph{Proof:}
As explained in the beginning of the section, we will couple between the measures $f$ and $g$ in means of coupling between the processes $\Gamma_t(f)$ and $\Gamma_t(g)$. To that end, we define, as in (\ref{contloc}),
\begin{equation}
F_0(x) = 1, ~~~ d F_t(x) = \langle A_t^{-1/2} d W_t, x - a_t \rangle F_t(x)
\end{equation}
where,
$$
a_t = \frac{\int_{\RR^n} x f(x) F_t(x) dx  }{ \int_{\RR^n} f(x) F_t(x) dx}
$$
is the barycenter of $f F_t$, and,
$$
A_t = \int_{\RR^n} (x - a_t) \otimes (x - a_t) f(x) F_t(x) dx
$$
is the covariance matrix of $f F_t$. As usual denote $f_t = F_t f$. \\
Next, we define,
$$
G_0(x) = 1, ~~~ d G_t(x) = \langle A_t^{-1/2} d W_t, x - b_t \rangle G_t(x)
$$
where,
$$
b_t = \frac{\int_{\RR^n} x g(x) G_t(x) dx  }{ \int_{\RR^n} g(x) G_t(x) dx},
$$
and denote $g_t(x) = g(x) G_t(x)$. \\ \\
Finally, we "interpolate" between the two processes by defining,
$$
H_0(x) = 1, ~~~ d H_t(x) = \langle A_t^{-1/2} d W_t, x - (a_t + b_t) / 2 \rangle,
$$
and,
$$
h_t(x) = H_t(x) H(f,g)(x).
$$
By a similar calculation to the one carried out in lemma \ref{basic1}, we learn that for all $t \geq 0$, $\int f_t(x) dx = \int g_t(x) dx = 1$.
Fix $x,y \in \RR^n$. An application of It\^{o}'s formula yields
$$
d \log f_t(x+y) = \langle x + y - a_t, A_t^{-1/2} d W_t \rangle - \frac{1}{2} |A_t^{-1/2}(x + y - a_t)|^2 dt,
$$
$$
d \log g_t(x-y) = \langle x - y - b_t, A_t^{-1/2} d W_t \rangle - \frac{1}{2} |A_t^{-1/2}(x - y - b_t)|^2 dt,
$$
and
$$
d \log h_t(x) = \left \langle x - \frac{a_t + b_t}{2}, A_t^{-1/2} d W_t \right \rangle - \frac{1}{2} |A_t^{-1/2} (x - (a_t + b_t) / 2) |^2 dt.
$$
Consequently,
$$
2 d \log h_t(x) \geq d \log f_t(x + y) + d \log g_t(x - y).
$$
It follows that,
$$
h_t(x) \geq H(f_t, g_t)(x).
$$
Define $S_t = \int_{\RR^n} h_t(x) dx$. The definition of $H_t$ suggests that $S_t$ is a martingale. By the Dambis / Dubins-Schwarz theorem, there exists a non-decreasing function $A(t)$ such that,
$$
S_t = K(f,g) + \tilde W_{A(t)}
$$
where $\tilde W_t$ is distributed as a standard Wiener process. Since $S_t \geq 1$ almost surely, it follows from the Doob's maximal inequality theorem that,
\begin{equation} \label{maxksmall}
\PP (G_t) \geq 1 - \epsilon / 2, ~~~ \forall s > 0.
\end{equation}
where,
\begin{equation} \label{defeventf}
G_t = \left \{ \max_{s \in [0, t] } S_s \leq \frac{2 K(f,g)}{\epsilon} \right \}.
\end{equation}
Next, define,
$$
F_t := \left \{ ||A_s||_{OP} < C K_n^2 (\log n) e^{-t}, ~~ \forall 0 \leq s \leq t \right \}.
$$
where $C$ is the same constant as in (\ref{deff}). Finally, denote $E_t = G_t \cap F_t$. By proposition 3.1 and equation (\ref{maxksmall}), $P(E_t) > 1 - \epsilon$ for all $t > 0$. Define a stopping time by the equation,
$$
\rho = \sup \{t | ~ E_t \mbox{ holds} \}.
$$
\medskip
Our next objective is to define the densities $\tilde f, \tilde g$ by, in some sense, neglecting the cases where $E_t$ does not hold.
We begin by defining the density $\tilde f_t$ by the following equation,
$$
\int_B \tilde f_t(x) dx = \EE \left [\mathbf{1}_{E_t} \int_B f_t(x) dx \right ],
$$
for all measurable $B \subset \RR^n$. Likewise, we define
$$
\int_B \tilde g_t(x) dx = \EE \left [\mathbf{1}_{E_t} \int_B g_t(x) dx \right ].
$$
Recall that $f(x) = \EE[f_t(x)]$ for all $x \in \RR^n$ and $t > 0$. It follows that,
$$
\int_{\RR^n} \tilde f_t(x) dx = \int_{\RR^n} \tilde g_t(x) dx = P(E_t) \geq 1 - \epsilon,
$$
and that 
$$
\tilde f_t(x) \leq f(x), ~~ \tilde g_t(x) \leq g(x), ~~ \forall x \in \RR^n.
$$
$$
~
$$
We construct a coupling between $\tilde f_t$ and $\tilde g_t$ by defining a measure $\mu_t$ on $\RR^n \times \RR^n$ using the formula
$$
\mu_t(A \times B) = \EE \left [ \mathbf{1}_{E_t} \int_{A \times B} f_t(x) g_t(y) dx dy \right ],
$$
for any measurable sets $A,B \subset \RR^n$. It is easy to check that $\tilde f_t$ and $\tilde g_t$ are the densities of the marginals of $\mu_t$ 
onto its first and last $n$ coordinates respectively. Thus, by definition of the Wasserstein distance,
$$
W_2(\tilde f_t, \tilde g_t) \leq \left ( \int_{\RR^n \times \RR^n} |x-y|^2 d \mu_t(x,y) \right )^{1/2} = 
$$
$$
\left ( \EE \left [\mathbf{1}_{E_t} \int_{\RR^n \times \RR^n} |x-y|^2 f_t(x) g_t(y) dx dy  \right ] \right )^{1/2} \leq
$$
$$
\left ( \EE \left [\mathbf{1}_{E_t} \left ( W_2(f_t, a_t) + W_2(g_t, b_t) + |a_t - b_t| \right )^2  \right ] \right )^{1/2}.
$$
Now, thanks to formula (\ref{atconverges}), we can take $T$ large enough (and deterministic) such that,
\begin{equation} \label{w2bytau}
W_2(\tilde f_T, \tilde g_T) \leq 2 \left ( \EE \left [\mathbf{1}_{E_T}  |a_T - b_T|^2  \right ] \right )^{1/2} + 1 \leq
\end{equation}
$$
2 \left ( \EE \left [|a_{T \wedge \rho} - b_{T \wedge \rho} |^2  \right ] \right )^{1/2} + 1.
$$
We will define $\tilde f := \tilde f_T$ and $\tilde g := \tilde g_T$. In view of the last equation, our main goal will be to attain a bound for the process $|a_t - b_t|$. A similar calculation to the one carried out in  (\ref{pathbc}) gives,
\begin{equation} \label{datdbt}
d a_t = A_t^{1/2} d W_t, ~~ d b_t = C_t A_t^{-1/2} d W_t.
\end{equation}
where,
$$
C_t = \int_{\RR^n} (x - b_t) \otimes (x - b_t) g_t(x) dx
$$
is the covariance matrix of $g_t$. Therefore,
$$
d |a_t - b_t|^2 = 2 \langle a_t - b_t, d a_t \rangle - 2 \langle a_t - b_t, d b_t \rangle +
$$
$$
\langle d a_t, d a_t \rangle + \langle d b_t, d b_t \rangle - 2 \langle d a_t, d b_t \rangle.
$$
The first two terms are martingale. We use the unique decomposition 
$$|a_t - b_t|^2 = M_t + N_t$$ 
where 
$M_t$ is a local martingale and $N_t$ is an adapted process of locally bounded variation.
We get,
$$
\frac{d}{dt} N_t = \langle d a_t - d b_t, d a_t - d b_t \rangle = 
$$
$$
\langle (A_t - C_t) A_t^{-1/2} d W_t, ( A_t - C_t) A_t^{-1/2} d W_t \rangle =
$$
$$
||A_t^{1/2} (I - A_t^{-1/2} C_t A_t^{-1/2})||_{HS}^2.
$$
By the Optional Stopping Theorem,
\begin{equation} \label{defdt}
\EE \left [ |a_{t \wedge \rho} - b_{t \wedge \rho}|^2 \right ] = \EE[N_{t \wedge \rho}] = \EE \left [\int_0^{t \wedge \rho} ||D_s||_{HS}^2 ds \right ]
\end{equation}
where $D_t = A_t^{1/2} (I - A_t^{-1/2} C_t A_t^{-1/2})$. Our next task is to use lemma \ref{transportation} to bound $||D_t||_{HS}$ under the assumption f$t < \tau$. \\ \\
We start by denoting the eigenvalues of the matrix $I - A_t^{-1/2} C_t A_t^{-1/2}$ by $\delta_i$, in decreasing order, and the eigenvalues of the matrix $A_t$ by $\lambda_i$, also in decreasing order.
By theorem 1 in \cite{tam}, 
\begin{equation} \label{traceineq}
||D_t||_{HS}^2 \leq \sum_{j=1}^n \lambda_j \delta_j^2.
\end{equation}
By lemma \ref{transportation}, we learn that
\begin{equation} \label{deltasmall}
\delta_j \leq \frac{C K(f_t,g_t)^{5} \tau_n j^{\kappa}}{\sqrt{j}}.
\end{equation}
Plugging this into (\ref{traceineq}) yields,
$$
||D_t||_{HS}^2 \leq C K(f_t,g_t)^{10} \tau_n^2 \sum_{j=1}^n \lambda_j j^{2 \kappa - 1}.
$$
Fix some constant $(1 - 2 \kappa) < \alpha < 1$, whose value will be chosen later. For now, we assume that $\kappa > 0$. Using H\"{o}lder's inequality, we calculate,
\begin{equation} \label{calcholder}
||D_t||_{HS}^2 \leq C K(f_t,g_t)^{10} \tau_n^2 \left ( \sum_{j=1}^n \lambda_j^{1 / (1-\alpha)} \right )^{1 - \alpha} \left (\sum_{j=1}^n j^{(2 \kappa - 1) / \alpha} \right)^{\alpha} \leq
\end{equation}
$$
C K(f_t,g_t)^{10} \tau_n^2 \left ( \lambda_1^{1 / (1-\alpha) - 1} \sum_{j=1}^n \lambda_j \right )^{1 - \alpha} \left (1 + \int_1^n t^{(2 \kappa - 1) / \alpha} \right )^{\alpha} \leq
$$
$$
C K(f_t,g_t)^{10} \tau_n^2 \lambda_1^{\alpha} (\beta n)^{1 - \alpha}  \left (n^{(2 \kappa - 1) / \alpha + 1} + 2 \right )^\alpha \left ( \frac{1}{(2 \kappa - 1) / \alpha + 1} \right )^\alpha
$$
where $\beta = \frac{1}{n} \sum_{j=1}^n \lambda_j$. Recall that $\alpha > (1 - 2 \kappa)$, which gives,
\begin{equation} \label{calcalpha}
\left ( n^{(2 \kappa - 1) / \alpha + 1} + 2 \right )^\alpha \leq 3 n^{\alpha} n^{2 \kappa - 1}.
\end{equation}
Take $\alpha$ such that $\epsilon = \alpha - (1 - 2 \kappa)$. Equations (\ref{calcholder}) and (\ref{calcalpha}) give,
$$
||D_t||_{HS}^2 \leq \frac{C'}{\epsilon}  K(f_t,g_t)^{10} \tau_n^2 \beta^{1 - \alpha} \lambda_1^{\alpha} n^{2 \kappa} \leq 
$$
$$
\frac{C''}{\epsilon}  K(f_t,g_t)^{10} \tau_n^2 \max(\beta, 1) \lambda_1^{1 - 2 \kappa + \epsilon} n^{2 \kappa}.
$$
Recall that we assume that $t < \tau$. By the definition of $\tau$, we
get $\lambda_1 \leq C \tau_n^2 n^{2 \kappa} \log n$ and $K(f_t, g_t) \leq 2 K(f,g) / \epsilon$.
Part (ii) of proposition \ref{mainsec3} implies $\EE[\beta] \leq 1$. 
Plugging these facts into the last equation gives,
$$
\EE \left [ ||D_t||_{HS}^2 ~ \right ] \leq \frac{C}{\epsilon^{11}} K(f,g)^{10} \tau_n^2 \left(\tau_n^2 n^{2 \kappa} \log n\right )^{1 - 2 \kappa + \epsilon} n^{2 \kappa} e^{-t} \leq 
$$
$$
\leq \frac{C'}{\epsilon^{11}} K(f,g)^{10} \tau_n^2 n^{4 \kappa - 4 \kappa^2 + \epsilon} e^{-t}.
$$
Finally, using equations (\ref{w2bytau}) and (\ref{defdt}), we conclude,
\begin{equation} \label{distsmall}
W_2(\tilde f_T, \tilde g_T)^2 \leq \EE \left [ \int_0^{T \wedge \rho} ||D_s||_{HS}^2 ds \right ] \leq 
\end{equation}
$$
\frac{C}{\epsilon^{11}} K(f,g)^{10} \tau_n^2 n^{4 \kappa - 4 \kappa^2 + \epsilon}.
$$
The proof is complete. \qed \\
\begin{remark}
In the above lemma, if we replace the assumption that $f$ is isotropic by the assumption that $f,g$ are log-concave with respect to the Gaussian measure, then following the same lines of proof while using proposition \ref{bakryemery}, one may improve the bound (\ref{smallw2}) and get,
$$
W_2(\tilde f, \tilde g) \leq C(\epsilon) K(f,g) \sqrt{\log n}.
$$
\end{remark}
\bigskip
We move on to the proof of theorem \ref{mainthm2}. \\ \\
\emph{Proof of theorem \ref{mainthm2}:}
Let $K,T$ be convex bodies of volume $1$ such that the covariance matrix of
$K$ is $L_k^2 Id$. Fix $\epsilon > 0$. Define,
$$
f(x) = 1_{K / L_K}(x) L_K^n, ~~~ g(x) = 1_{T / L_K}(x) L_K^n,
$$
so both $f$ and $g$ are probability measures and $f$ is isotropic. We have,
$$
K(f,g) = Vol_n \left (\frac{K+T}{2} \right ) = V.
$$
We use lemma \ref{mainlemmasec5}, which asserts that there exist two measures $\tilde f$, $\tilde g$, such that,
\begin{equation} \label{fgprop1}
\tilde f(x) \leq f(x), ~~ \tilde g(x) \leq g(x), ~~ \forall x \in \RR^n,
\end{equation}
\begin{equation} \label{fgprop2}
\int \tilde f(x) dx = \int \tilde g(x) dx \geq 1 - \epsilon
\end{equation}
and such that,
$$
W_2 (\tilde f, \tilde g) \leq \Theta
$$
where $\Theta = C(\epsilon) V^{5} \tau_n n^{2 (\kappa - \kappa^2) + \epsilon}$.
Since $\tilde g$ is supported on $T$, it follows that,
$$
\int_{K} d^2 (x, T / L_T) \tilde f(x) dx \leq \Theta^2
$$
where $d(x,T / L_T) = \inf_{y \in (T / L_T)} |x - y|$. Denote,
$$
K_\alpha = \{x \in K / L_K; ~ d(x,T) \geq \alpha \Theta \}.
$$ 
It follows from  Markov's inequality and from (\ref{fgprop1}) and (\ref{fgprop2}) that,
$$
Vol_n(K_\alpha) \leq L_K^{-n} \left (\epsilon + \frac{1}{\alpha^2} \right ).
$$
Finally, taking $\delta = L_K \Theta / \sqrt{\epsilon}$ gives
\begin{equation} \label{markov}
Vol_n(K \setminus T_{\delta}) \leq 2 \epsilon.
\end{equation}
This completes the proof. \qed \\ \\
\section{Tying up loose ends}

We begin the section with the proof of lemma \ref{ksigma} which gives an upper bound for the constant $K_n$ in terms of $\tau_n$ and $\kappa$. \\ \\
\emph{Proof of lemma \ref{ksigma}:}
Let $X$ be an isotropic, log concave random vector in $\RR^n$, and fix
$\theta \in \Sph$.
Denote $A = \EE[X \otimes X \langle X, \theta \rangle]$. Our goal is to show,
$$
||A||_{HS}^2 \leq C \sum_{k=1}^n \frac{\sigma_k^2}{k}.
$$
Let $k \leq n$ and let $E_k$ be a subspace of dimension $k$.
Denote $P(X) = Proj_{E_k}(X)$ and $Y = |P(X)| - \sqrt k$. By definition of $\sigma_k$,
$$
Var[Y] \leq \sigma_k^2
$$
Note that, by the isotropicity of $X$, $\EE[|P(X)|^2] = k$. It easily follows that,
$$
Var[|P(X)|^2] \leq C k Var[Y] \leq C k \sigma_k^2.
$$
Using the last inequality and applying Cauchy-Schwartz gives,
$$
\left | \EE [\langle X, \theta \rangle |P(X)|^2] \right | \leq \sqrt{Var[\langle X, \theta \rangle] Var[|P(X)|^2]} \leq C \sqrt{k} \sigma_k
$$
or, in other words,
$$
\left | Tr[Proj_{E_k} A Proj_{E_k}] \right | \leq C \sqrt k \sigma_k.
$$
Let $\lambda_1,...,\lambda_\ell$ be the non-negative eigenvalues of $A$ in decreasing order. The last inequality implies that the matrix $Proj_{E_k} A Proj_{E_k}$ has at least one eigenvalue smaller than $C \sqrt{\frac{1}{k}} \sigma_k$. Consequently, by taking $E_k$ to be the subspace spanned by the $k$ first corresponding eigenvectors we learn that
$$
\lambda_k^2 < C \frac{\sigma_k^2}{k}, ~~ \forall k \leq \ell.
$$
In the same manner, if $\zeta_1,...,\zeta_{n-\ell}$ are the negative eigenvalues of $A$, one has $\zeta_k^2 < C \frac{\sigma_k^2}{k}$. \\
We can thus calculate,
$$
||A||_{HS}^2 = \sum_{k=1}^\ell \lambda_k^2 + \sum_{k=1}^{n - \ell} \zeta_k^2 \leq 2 C \sum_{k=1}^n \frac{\sigma_k^2}{k}.
$$
The proof is complete. \qed \\ \\

Next, in order to provide the reader with a better understanding of the constant $K_n$, we introduce two new constants. 
First, define
$$
Q_n^2 = \sup_{X,Q} \frac{Var[Q(X)]}{\EE \left [|\nabla Q(X)|^2 \right ]}
$$
where the supremum runs over all isotropic log-concave random vectors, $X$, and all quadratic forms $Q(x)$. Next, define
$$
R_n^{-1} = \inf_{\mu, E} \frac{\mu^+(E)}{\mu(E)}
$$
where $\mu$ runs over all isotropic log-concave measures and $E$ runs over all ellipsoids with $\mu(E) \leq 1/2$. \\
\begin{fact}
There exist universal constants $C_1, C_2$ such that 
$$K_n \leq C_1 Q_n \leq C_2 R_n.$$
\end{fact}
The proof of the right inequality is standard and uses the coarea formula and the Cauchy-Schwartz inequality. We will prove the left inequality. To that end, fix an isotropic log-concave random vector $X$, denote $A = \EE[X \otimes X X_1]$. We have,
$$
||A||_{HS} = \sup_{B} \frac{Tr(BA)}{||B||_{HS}}
$$
where $B$ runs over all symmetric matrices. Let $B$ be a symmetric matrix. Fix coordinates under which $B$ is diagonal, and write $X = (X_1,...,X_n)$ 
and $B = diag \{a_1,..,a_n\}$. Define $Q(x) = \langle B x, x \rangle$. We have,
$$
Tr(BA) = \EE \left [X_1 \sum_{i=1}^n a_i X_i^2 \right ] \leq \sqrt{\EE \left [X_1^2 \right ]} \sqrt{Var \left [\sum_{i=1}^n a_i X_i^2 \right ]} =
$$
$$
\sqrt{Var[Q(X)]} \leq  \sqrt{2 Q_n^2 \sum_{i=1}^n a_i^2 \EE[X_i^2 ]} = \sqrt{2} Q_n ||B||_{HS}.
$$
So,
$$
||A||_{HS} \leq \sqrt{2} Q_n.
$$
This shows that $K_n \leq C Q_n$.
\begin{remark}
We suspect that there exists a universal constant $C>0$ such that $K_n \leq C \sigma_n$, but we are unable to prove that assertion.
\end{remark}
\bigskip
We move on to the proof of lemma \ref{basic1.5}. \\
\emph{Proof of lemma \ref{basic1.5}:} \\
Throughout the proof, all the constants $c, c_1, c_2,...$ may depend only on the dimension $n$. 
Recall that $f(x)$ is assumed to be isotropic and log-concave. It is well-known that there exist two constants $c_1, c_2>0$, such that 
$$
f(|x|) \geq c_1, ~~ \forall |x| \leq c_2.
$$
(see for example \cite[Theorem 5.14]{LV}). Define $g(x) = c_1 \mathbf{1}_{\{ |x| \leq c_2 \} }$. 
It is also well-known (see for example \cite[Lemma 5.7]{LV}) that there exist two constants $c_3, c_4 > 0$ such that
$$
\int_{\RR^n} f(x) e^{\langle x, y \rangle } dx \leq c_3, ~~ \forall |y| < c_4,
$$
which implies that whenever $|c| < c_4$ and $B$ is positive semi-definite,
$$
V_f(c, B) = \int_{\RR^n} e^{\langle c, x \rangle - \frac{1}{2} \langle B x, x \rangle } f(x) dx \leq c_3.
$$
It follows that for all $|c| < c_4$ and $B \leq Id$ (in the sense of positive matrices), one has
\begin{equation} \label{covmatbig}
A_f(c,B) \geq 
\end{equation}
$$
c_3^{-1} \int_{\RR^n} (x - a_f(c,B)) \otimes (x - a_f(c,B)) e^{- c_4 |x| - \frac{1}{2} |x|^2 } g(x) dx \geq 
$$
$$
c_3^{-1} c_1 \int_{\{|x| \leq c_2\}} x \otimes x e^{- c_4 |x| - \frac{1}{2} |x|^2 } dx = c_5 Id
$$
for some constant $c_5>0$.
Define the stopping times,
$$
T_1 = \sup \{t > 0; ~|c_t| < c_4  \}, ~~ T_2 = \sup \{t > 0; ~B_t \geq Id  \}, ~~ T = \min(T_1,T_2).
$$
Note that according to (\ref{covmatbig}),
\begin{equation} \label{covmatbig2}
A_t \geq c_5 Id, ~~ \forall t \leq T,
\end{equation}
so the lemma would be concluded if we manage to show that
\begin{equation} \label{Tlarge}
\PP(T > c) > c
\end{equation}
for some constant $c>0$. \\ \\
Define the event $E = \{T_2 \leq T_1\}$. Whenever $E$ holds, we have
the following: First, using (\ref{covmatbig}), 
$$
A_t \geq c_5 Id, ~~ \forall t \leq T_2.
$$
Recall that $\frac{d}{dt} B_t = A_t^{-1}$. It follows that 
$$
B_t \leq c_5^{-1} t,  ~~ \forall t \leq T_2.
$$
By taking $t = T_2$ in the last equation, we learn that $T = T_2 \geq c_5$ whenever $E$ holds, so
$$
T_2 \leq T_1 \Rightarrow T \geq c_5.
$$
Therefore, it is enough to prove that $P(T_1 > c) > c$ for some $c>0$. Furthermore, in the following we 
are able to assume that $P(E) \leq 0.1$. \\ \\
To that end, consider the defining equation (\ref{stochastic1}) and use It\^{o}'s formula to attain
\begin{equation} \label{dct2}
d |c_t|^2 = 2 \langle c_t, A_t^{-1/2} dW_t \rangle + 2 \langle A_t^{-1} a_t, c_t \rangle dt + ||A_t^{-1/2}||_{HS}^2 dt.
\end{equation}
Define the process $e_t$ by the equations,
$$
e_t = 0, ~~ d e_t = 2 \langle c_t, A_t^{-1/2} dW_t \rangle.
$$
Using (\ref{covmatbig2}), we deduce that whenever $t < T$, one has
\begin{equation} \label{qvbound}
[e]_t = 4 \int_0^t \langle A_t^{-1/2} c_t, A_t^{-1/2} c_t \rangle \leq 4 c_4^2 c_5^{-1} t.
\end{equation}
Using the Dambis / Dubins-Schwartz theorem, we know that there exists a standard Wiener process $\tilde W_t$
such that $e_t$ has the same distribution as $\tilde W_{[e]_t}$. An elementary property of the standard Wiener process
is that there exists a constant $c_6 > 0$, such that 
\begin{equation}
\PP(F) \geq 0.9
\end{equation}
where 
$$
F = \left \{ \max_{0 \leq s \leq c_6} \tilde W_s \leq c_4^2 / 2  \right \}.
$$
Define $\delta = \min \left (T, \frac{c_6}{4 c_4^2 c_5^{-1}} \right)$. Note that, by (\ref{qvbound}),
\begin{equation} \label{fsubset}
F \subseteq \left \{ \max_{0 \leq t \leq \delta}  e_t \leq c_4^2 / 2 \right \}.
\end{equation}
Another application of (\ref{covmatbig2}), this time with the assumption (\ref{compactsupport}) gives,
\begin{equation} \label{predbound}
\int_0^{t} ||A_s^{-1/2}||_{HS}^2 ds + 2 \left | \langle A_s^{-1} a_s, c_s \rangle \right | ds \leq n c_5^{-1} (1 + c_4) t \leq c_7 t
\end{equation}
for all $t < T$, and for a constant $c_7$. By plugging (\ref{fsubset}) and (\ref{predbound}) into (\ref{dct2}), we learn that whenever $F$ holds, 
one has
$$ 
|c_t|^2 \leq c_4^2 / 2 + c_7 t, ~~ \forall t \leq \delta.
$$
If we assume that $\delta = T_1$, the above gives $c_4^2 = |c_{T_1}|^2 \leq c_4^2 / 2 + c_7 {T_1}$ which implies $T \geq \frac{c_4^2}{2 c_7}$ (here, we used the assumption that $T_1 \leq T_2$). Thus, whenever $F \cap E^C$ holds, we have $T = T_1 \geq \min \left (\frac{c_4^2}{2 c_7}, \frac{c_6}{4 c_4^2 c_5^{-1}} \right )$. A union bound gives $\P(F \cap E^C) \geq 0.8$. The lemma is complete. \qed \\ \\

\bigskip

{\small \noindent \it e-mail address: roneneldan@gmail.com}

\vfill \hfill \today


\begin{thebibliography}{GGM}


\bibitem[ABP]{ABP} M. Anttila, K. Ball, I. Perissinaki, {\it The
central limit problem for convex bodies}. Trans. Amer. Math. Soc.,
355, no. 12, (2003), 4723--4735.

\bibitem[BE]{BE} D. Bakry and M. Emery, {\it Diffusions hypercontractives}, in S\'{e}minaire de probabilit\'{e}s, XIX,
1983/84, vol. 1123 of Lecture Notes in Math., Springer, Berlin, 1985, pp. 177–206.

\bibitem[BN]{bn}
K. Ball V.H. Ngyuen, {\it Entropy jumps for random vectors with log-concave
density and spectral gap.} Preprint.


\bibitem[Bo]{Bobkov} S. Bobkov, {\it On isoperimetric constants for log-concave probability distributions,
in Geometric Aspects of Functional Analysis } Israel Seminar 2004-2005, Springer
Lecture Notes in Math. 1910 (2007), 8188.

\bibitem[BK]{BK} S. Bobkov, A. Koldobsky, {\it
On the central limit property of convex bodies. } Geometric aspects
of functional analysis, Lecture Notes in Math., 1807, Springer,
Berlin, (2003), 44--52.

\bibitem[Bou]{bou3} Bourgain, J.,
{\it On the distribution of polynomials on high-dimensional convex
sets. } Geometric aspects of functional analysis, Israel seminar
(1989--90), Lecture Notes in Math., 1469, Springer, Berlin, (1991),
127--137.

\bibitem[Dis]{diskant} Diskant, V. I., {\it Stability of the Solution of the Minkowski
Equation} (in Russian). Sibirsk. Mat. 14 (1973), 669--673, 696.
English translation in Siberian Math. J. 14 (1973), 466-469.

\bibitem[Dur]{durett} R Durrett, {\it Stochastic Calculus: A Practical Introduction} Cambdidge university press, 2003.

\bibitem[EK1]{EK1} Eldan, R., Klartag, B., {\it Approximately gaussian marginals and the hyperplane
conjecture}. Proc. of a workshop on ``Concentration, Functional
Inequalities and Isoperimetry'', Contemporary Math., vol. 545, Amer.
Math. Soc., (2011), 55--68.

\bibitem[EK2]{EK2} Eldan, R., Klartag, B., {\it Dimensionality and the stability of the
Brunn-Minkowski inequality}. Annali SNS, 2011.

\bibitem[Fl1]{fleury} B. Fleury, {\it Concentration in a thin euclidean shell for log-concave measures }, J.
Func. Anal. 259 (2010), 832841.

\bibitem[Fl2]{fleury2} B. Fleury, {\it Poincar\'{e} inequality in mean value for Gaussian polytopes}, 
Probability theory and related fields, Volume 152, Numbers 1-2, 141-178.

\bibitem[FMP1]{FMP1} Figalli, A., Maggi, F., Pratelli, A., {\it A refined Brunn-Minkowski
inequality for convex sets.} Ann. Inst. H. Poincar\'e Anal. Non
Lin\'eaire,  vol. 26, no. 6, (2009), 2511–-2519.

\bibitem[FMP2]{FMP2} Figalli, A., Maggi, F., Pratelli, A., {\it A mass transportation approach
to quantitative isoperimetric inequalities.} Invent. Math., vol.
182,  no. 1, (2010), 167–- 211.

\bibitem[Gu-M]{GM} O. Guedon, E. Milman, {\it Interpolating thin-shell and sharp large-deviation estimates for isotropic log-concave measures}, 2010

\bibitem[Gr-M]{GrMil} M. Gromov and V. D. Milman. {\it A topological application of the isoperimetric inequality.}
Amer. J. Math., 105(4):843–854, 1983

\bibitem[Groe]{groemer} Groemer, H., {\it On the Brunn–-Minkowski theorem}. Geom. Dedicata, vol. 27, no. 3, (1988), 357–-371.

\bibitem[Gros1]{Gross} L. Gross { \it Logarithmic Sobolev inequalities} Amer. J. Math. 97 (1975), no. 4, 1061-1083

\bibitem[Gros2]{Grossbook} L. Gross, { \it Logarithmic Sobolev inequalities and contractivity properties of semigroups, Dirichlet forms} Varenna, 1992, 54-88, Lecture Notes in Math., 1563, Springer, Berlin, 1993.

\bibitem[Leh]{Lehec} J. Lehec, {\it Representation formula for the entropy and functional inequalities.} arXiv: 1006.3028, 2010.

\bibitem[KX]{KX} G. Kallianpur, J. Xiong, {\it Stochastic Differential Equations in Infinite Dimensional Spaces}
Institute of mathematical statistics, lecture notes - monograph series. California, USA, 1995.

\bibitem[K1]{K1} Klartag, B., {\it
A central limit theorem for convex sets.} Invent. Math., 168,
(2007), 91--131.

\bibitem[K2]{K2} Klartag, B., {\it
Power-law estimates for the central limit theorem for convex sets.}
J. Funct. Anal., Vol. 245, (2007), 284--310.

\bibitem[K3]{K_uncond} Klartag, B., {\it A Berry-Esseen type inequality for convex bodies with an
unconditional basis. } Probab. Theory Related Fields, vol. 145, no.
1-2, (2009), 1–-33.

\bibitem[K4]{K4} Klartag, B., {\it Power-law estimates for the central limit theorem for convex sets. } wherever, 2007.

\bibitem[KLS]{KLS}
R. Kannan, L. Lov\'{a}sz, and M. Simonovits. Isoperimetric problems for convex bodies and a localization lemma. Discrete Comput. Geom., 13(3-4):541–559, 1995

\bibitem[L]{ledoux}
M. Ledoux, {\it Spectral gap, logarithmic Sobolev constant, and geometric bounds.}
Surveys in differential geometry. Vol. IX, 219-240, Surv. Differ. Geom., IX, Int. Press, Somerville, MA, 2004.

\bibitem[LV]{LV} L. Lov\'asz and S. Vempala, {\it The geometry of logconcave functions and sampling algorithms.}  Random Structures \& Algorithms, Vol. 30, no. 3, (2007),  307--358.

\bibitem[Mil1]{Mil}
E. Milman, {\it On the role of Convexity in Isoperimetry, Spectral-Gap and Concentration},  Invent. Math. 177 (1), 1-43, 2009.

\bibitem[Mil2]{Mil2}
E. Milman, {\it Isoperimetric Bounds on Convex Manifolds},  Contemporary Math., proceedings of the Workshop on "Concentration,Functional Inequalities and Isoperimetry" in Florida, November 2009.

\bibitem[Ok]{oksendal} B. Oksendal { \it Stochastic Differential Equations: An Introduction with Applications.} Berlin: Springer. ISBN 3-540-04758-1, (2003).

\bibitem[Oss]{osserman} R. Osserman, {\it Bonnesen-style isoperimetric inequalities. } Amer.
Math. Monthly, 86, no. 1,  (1979), 1–-29.

\bibitem[Pis]{pisier} G. Pisier, {\it The volume of convex bodies and Banach space geometry. }
Cambridge Tracts in Mathematics, 94. Cambridge University Press,
Cambridge, 1989.

\bibitem[Seg]{segal} A. Segal, {\it Remark on Stability of Brunn-Minkowski and Isoperimetric
Inequalities for Convex Bodies}. To appear in Gafa Seminar notes.

\bibitem[Sud]{sudakov}
V.N. Sudakov, {\it Typical distributions of linear functionals in finite-dimensional
spaces of high dimension.} (Russian) Dokl. Akad. Nauk SSSR 243 (1978), no. 6,
1402-1405.

\bibitem[T]{tam} T. Tam, {\it On Lei-Miranda-Thompson's result on singular values and diagonal elements}.
Linear Algebra and Its Applications, 272 (1998), 91-101.

\bibitem[Vil]{villani} Villani, C., {\it Topics in optimal transportation. }
Graduate Studies in Mathematics, 58. American Mathematical Society, Providence, RI, 2003.

\end{thebibliography}
\end{document}